\documentclass[12pt]{article}
\usepackage{e-jc}

\usepackage{amsmath,tikz}
\usepackage[mathscr]{euscript}

\theoremstyle{definition}
\newtheorem*{sub*}{Subdivision recurrence}
\newtheorem*{tri*}{Triangle recurrence}
\newtheorem*{defn*}{Definition}

\DeclareMathOperator{\aff}{aff}
\DeclareMathOperator{\NVol}{NVol}
\DeclareMathOperator{\vol}{vol}
\DeclareMathOperator{\conv}{conv}

\DeclareMathOperator{\PQ}{PQ}

\newcommand{\alphatri}{\alpha^{\triangle}}
\newcommand{\betatri}{\beta^{\triangle}}
\newcommand{\gammatri}{\gamma^{\triangle}}

\newcommand{\scrA}{\mathscr{A}}
\newcommand{\scrAtri}{\mathscr{A}^{\triangle}}
\newcommand{\scrB}{\mathscr{B}}
\newcommand{\scrBtri}{\mathscr{B}^{\triangle}}
\newcommand{\scrC}{\mathscr{C}}
\newcommand{\scrCtri}{\mathscr{C}^{\triangle}}
\newcommand{\scrF}{\mathscr{F}}

\newcommand{\calN}{\mathcal{N}}
\newcommand{\fkD}{\mathfrak{D}}

\newcommand{\RR}{\mathbb{R}}
\newcommand{\ZZ}{\mathbb{Z}}

\newcommand{\col}{\mathbin{:}}

\newcommand{\wkd}[1]{{#1}^{(*)}}
\newcommand{\ewd}[1]{{#1}^{(\ast\ast)}}
\newcommand{\APQ}{\nabla^{\PQ}}

\dateline{August 7, 2020}{March 11, 2022}{TBD}
\Copyright{The authors. Released under the CC BY-ND license (International 4.0).}

\title{Computing Volumes of Adjacency Polytopes via Draconian Sequences}

\author{Robert Davis\thanks{Supported by NSF grant DMS-1922998.}\\
\small Department of Mathematics\\[-0.8ex]
\small Colgate University\\[-0.8ex] 
\small Hamilton, New York, U.S.A.\\
\small\tt rdavis@colgate.edu\\
\and
Tianran Chen\thanks{Supported by NSF grant DMS-1923099.}\\
\small Department of Mathematics\\[-0.8ex]
\small Auburn University -- Montgomery\\[-0.8ex]
\small Montgomery, Alabama, U.S.A.\\
\small\tt tchen1@aum.edu}

\begin{document}

\maketitle 

\begin{abstract}
Adjacency polytopes appear naturally in the study of nonlinear emergent phenomena in complex networks. The ``PQ-type'' adjacency polytope, denoted $\nabla^{\mathrm{PQ}}_G$ and which is the focus of this work, encodes rich combinatorial information about power-flow solutions in sparse power networks that are studied in electric engineering. Of particular importance is the normalized volume of such an adjacency polytope, which provides an upper bound on the number of distinct power-flow solutions.
        
In this article we show that the problem of computing normalized volumes for $\nabla^{\mathrm{PQ}}_G$ can be rephrased as counting $D(G)$-draconian sequences where $D(G)$ is a certain bipartite graph associated to the network. We prove recurrences for all networks with connectivity at most $1$ and, for $2$-connected graphs under certain restrictions, we give recurrences for subdividing an edge and taking the join of an edge with a new vertex. Together, these recurrences imply a simple, non-recursive formula for the normalized volume of $\nabla^{\mathrm{PQ}}_G$ when $G$ is part of a large class of outerplanar graphs; we conjecture that the formula holds for all outerplanar graphs. Explicit formulas for several other (non-outerplanar) classes are given. Further, we identify several important classes of graphs $G$ which are planar but not outerplanar that are worth additional study.
\end{abstract}

\tableofcontents


\section{Introduction and background}

Let $G = (V(G),E(G))$ be a simple graph on $[N] = \{1,\dots,N\}$.
We use $e_1,\dots,e_N$ to denote the standard basis vectors of $\RR^N$.
The \emph{PQ-type adjacency polytope} of $G$ is defined to be
\[
    \APQ_G = \conv\{ (e_i,e_j) \in \RR^{2N} \mid ij \in E(G) \text{ or } i=j\}
\]
where $\conv(S)$ denotes the convex hull of elements of $S$.
Its \emph{normalized volume}, defined by $\NVol(\APQ_G) = \dim(\APQ_G)! \vol(\APQ_G)$ where $\vol(P)$ is the relative volume of $P$, is always a positive integer.

The study of PQ-type adjacency polytopes was introduced in \cite{ChenMehtaPQ}, 
motivated by the engineering problem known as \emph{power-flow study} (or \emph{load-flow study}).
This study models the balance of
electric power on a network of power generation or delivery ``buses''.
Of particular importance are the alternating current (AC) variations,
which produce nonlinear equations that are notoriously difficult to analyze.
In the AC model for a power network with buses labeled as $1,\dots,N$, 
the voltage on each bus is expressed as a complex variable
$v_i = x_i + \mathbf{i} y_i$ whose absolute value represents the voltage magnitude
and whose argument encodes the phase of the AC experienced on the bus.
The interaction among buses is modeled by a graph $G$ whose nodes 
represent the buses and whose edges represent the junctions.
Kirchhoff's circuit laws give rise to an idealized balancing condition
for the power injected, power generated, and power consumed
on each bus, which can be expressed as the system of nonlinear equations
\begin{equation}
    S_i = \sum_{j=1}^N \overline{Y}_{ij} v_i \overline{v}_j
    \quad\text{for } i = 2,\dots,N,
\end{equation}
where $S_i = P_i + \mathbf{i} Q_i$ is a complex representation of the 
real and reactive power,
$Y_{ij}$, known as \emph{nodal admittance}, describes the connection between
the $i$ and $j$ buses,
and $\overline{Y}_{ij}$ and $\overline{v}_j$ denote the complex conjugate
of $Y_{ij}$ and $v_j$ respectively.
By dropping the conjugate constraints between $v_i$ and $\overline{v}_i$,
we obtained the algebraic version of this system, 
known as the algebraic power-flow equations.
It was shown that the maximum number of nontrivial complex solutions
this system has is bounded by the normalized volume of $\APQ_G$.

We take care to call the adjacency polytopes within this paper PQ-type, since a related construction is sometimes called an adjacency polytope; see, for example, \cite{BraunBruegge, ChenDavisMehta, ManyFaces}.
This alternate construction, motivated by counting equilibrium solutions to a network of interconnected oscillators, relies on a particular change of variables that is not available here. 
In engineering terms, this alternate construction arises from PV-type buses.

In this article we show that the normalized volume of $\APQ_G$ can be described in terms of sequences of nonnegative integers related to the \emph{Dragon Marriage Problem}: a variant of Hall's Matching Theorem that has far-reaching applications and spawned the study of generalized permutohedra \cite{ PostnikovEtAl, Postnikov2009}.
We establish this relationship in Section~\ref{sec:translation} and show how it can be immediately exploited to compute normalized volumes of some PQ-type adjacency polytopes when $G$ is nontrivial.

We explore this connection more deeply in Section~\ref{sec: recurrences} where we establish several recurrences.
Namely, we provide recurrences for all graphs with connectivity at most $1$, that is, any graph that is disconnected or has a cut-vertex.
These directly imply a simple formula for $\NVol(\APQ_G)$ whenever $G$ is a forest. 

Sections~\ref{subsec: subdivision} and \ref{subsec: triangle} consider two operations on a graph: subdivision of an edge $e$ and replacing $e$ with the join of $e$ and a new vertex. 
Under certain conditions, these operations lead to the following two recurrences that are stated simply but nontrivial to prove.

\begin{sub*}[see Theorem~\ref{thm: subdiv recurrence}]
	Let $G$ be a $2$-connected graph and let $e=uv$ be an edge.
	Denote by $G \col e$ the graph obtained by subdividing $e$.
	If $\deg_G(u)=2$ and the neighbors of $u$ are neighbors of each other, then
	\[
		\NVol(\APQ_{G \col e}) = 2\NVol(\APQ_G) + \NVol(\APQ_{G \setminus e}).
	\]
\end{sub*}

\begin{tri*}[see Theorem~\ref{thm: triangle recurrence}]
	Let $G$ be any connected graph and let $e = uv$ be an edge with $\deg_{G}(u) = 2$.
	If $\deg_G(v) = 2$ or if the neighbors of $u$ are neighbors of each other, then
	\[
		\NVol(\APQ_{G \triangle e}) = 3\NVol(\APQ_G).
	\]
\end{tri*}

Section~\ref{sec: recurrences} concludes by applying the recurrences to establish a closed, non-recursive formula for $\NVol(\APQ_G)$ for a large class of outerplanar graphs; we conjecture that this formula holds for all outerplanar graphs.
The final section addresses several classes of graphs which are planar but not outerplanar.
First, we give results for a complete bipartite graph where one partite set has just two elements.
Then we consider the classes of \emph{wheel graphs} and \emph{series-parallel graphs}, which are natural points of further study and will likely require a refinement of the techniques within this article or alternate techniques altogether.


\section{Notation, background, and translating to draconian sequences}\label{sec:translation}

Before we prove our results, we will establish assorted notation that will be needed throughout this work.
Additional notation will be introduced as needed.
First, if $e$ is an edge of $G$ with endpoints $u$ and $v$, we will write $e = uv$ or $e = vu$ whenever possible.
When additional clarity is helpful we may alternately write $e = \{u,v\}$ or $e = \{v,u\}$.

If $X \subseteq V(G)$, then we use $G-X$ to denote the graph obtained from deleting the vertices of $X$ as well as any edge that is incident to some vertex in $X$. 
If $X = \{v\}$, then we will just write $G - v$.
Similarly, if $S$ is a set of edges, then we use $G\setminus S$ to denote the graph with the edges in $S$ deleted; if $S = \{e\}$, then we just write $G - e$.
If $X \subseteq V(G)$, then we use $G[X]$ to denote the subgraph of $G$ induced by $X$.
Lastly, if $H$ is a graph, then we use $G \vee H$ to denote the \emph{join} of $G$ and $H$, that is, the graph with vertex set $V(G) \cup V(H)$ and edge set
\[
	E(G) \cup E(H) \cup \{uv \mid u \in V(G), v \in V(H)\}.
\]
For a positive integers $M,N$, let $K_N$ denote the complete graph on $[N]$ and let $K_{M,\overline N}$ denote the complete bipartite graph with partite sets $[M]$ and $[\overline{N}] = \{\overline 1, \dots, \overline{N}\}$.
Let $\calN_G(v)$ denote the set of vertices of $G$ adjacent to $v$.
Keeping this notation in mind, we may now begin in earnest.

In \cite{Postnikov2009}, Postnikov investigated the \emph{Dragon Marriage Problem}, providing a generalization of Hall's Matching Theorem for bipartite graphs.
In the Dragon Marriage Problem, a small medieval village is home to $n$ grooms and $n+1$ brides, some pairs of whom would form compatible marriages. 
Suppose we know all pairs of compatible grooms and brides. 
One day, a dragon arrives in the village and kidnaps a bride.
What compatibility conditions among the original set of grooms and brides will guarantee that those who remain can still be entirely paired by compatible marriages?
In graph-theoretic terms, and more generally, consider an $X,Y$-bigraph $G$ such that $|Y| = |X|+1$.
What are necessary and sufficient conditions on $G$ so that $G-y$ has a perfect matching regardless of choice of $y \in Y$?
The answer relies on the following.

\begin{definition}\label{def: draconian sequence}
	Let $G \subseteq K_{N,\overline N}$.
	Call $(a_1,\dots,a_N) \in \ZZ_{\geq 0}^N$ a \emph{$G$-draconian sequence} if $\sum a_i = N-1$ and, for any $1 \leq i_1 < i_2 < \cdots < i_k \leq N$,
	\begin{equation}\label{eq: draconian inequality}
		a_{i_1} + \cdots + a_{i_k} < \left|\bigcup_{j=1}^k \calN_G(i_j)\right|.
	\end{equation}
	We will say that a sequence satisfying \eqref{eq: draconian inequality} satisfies the \emph{$G$-draconian inequality corresponding to $i_1,\dots,i_k$}.
\end{definition}

Postnikov proved \cite[Proposition 5.4 and Definition 9.2]{Postnikov2009} that a matching that covers $X$ exists exactly when a $G$-draconian sequence exists.
He then goes on to compute volumes of certain polyhedra as sums over the set of $G$-draconian sequences.
At the moment, it may be completely unclear how draconian sequences are useful to us; the rest of this section is dedicated to clarifying the connection.

\begin{definition}\label{def: root}
	Given a graph $G \subseteq K_{M,\overline N}$, let $Q_G$ denote the \emph{root polytope} 
	\[
		Q_G = \conv\{e_i - e_{\overline j} \mid \{i,\overline j\} \in E(G)\} \subseteq \RR^M \times \RR^{\overline N},
	\]
	where $\RR^{\overline N}$ denotes the real vector space with standard basis vectors $e_{\overline 1},\dots,e_{\overline N}$.
\end{definition}

It turns out that we can describe $\APQ_G$ as a root polytope for an appropriate choice of graph.

\begin{definition}
	Let $G$ be a simple graph on $[N]$.
	Define $D(G)$ to be the subgraph of $K_{N,\overline N}$ with edges $\{i, \overline i\}$ for each $i \in [N]$ and $\{i, \overline j\}$ and $\{j, \overline i\}$ for each edge $ij$ in $G$.
\end{definition}

\begin{figure}
\begin{center}
\begin{tikzpicture}
\begin{scope}[every node/.style={circle,fill,inner sep=0pt,minimum size=2mm}]
	\node[label=$1$] (A) at (-1,0) {};
	\node[label=$2$] (B) at (0.8,0) {};
	\node[label=right:$3$] (C) at (2,1) {};
	\node[label=right:$4$] (D) at (2,-1) {};
\end{scope}
\node (E) at (0,-1.5) {};

\draw[thick] (A) -- (B) -- (C) -- (D) -- (B);
\end{tikzpicture}
\hspace{1in}
\begin{tikzpicture}

\begin{scope}[every node/.style={circle,fill,inner sep=0pt,minimum size=2mm}]
	\node[label=left:$1$] (A) at (0,0) {};
	\node[label=left:$2$] (B) at (0,1) {};
	\node[label=left:$3$] (C) at (0,2) {};
	\node[label=left:$4$] (D) at (0,3) {};
	\node[label=right:$\overline{1}$] (AA) at (1,0) {};
	\node[label=right:$\overline{2}$] (BB) at (1,1) {};
	\node[label=right:$\overline{3}$] (CC) at (1,2) {};
	\node[label=right:$\overline{4}$] (DD) at (1,3) {};
\end{scope}

\draw[thick] (A) -- (AA) -- (B) -- (BB) -- (A);
\draw[thick] (B) -- (CC) -- (C) -- (BB) -- (D) -- (CC);
\draw[thick] (D) -- (DD) -- (C);
\draw[thick] (B) -- (DD);
\end{tikzpicture}
\end{center}
\caption{A graph $G$, left, and its corresponding bipartite graph $D(G)$, right.}\label{fig: D(G)}
\end{figure}
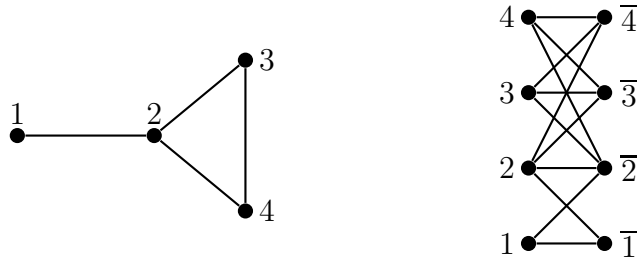

As an example, let $G$ be the graph on $[4]$ with edges $12$, $23$, $34$, $24$.
Then $D(G)$ is the bipartite graph with vertices $\{1,2,3,4, \overline 1, \overline 2, \overline 3, \overline 4\}$ and edges
$1\overline 1$, $1\overline 2$, $2\overline 1$, $2\overline 2$, $2 \overline 3$, $2 \overline 4$, $3 \overline 2$, $3 \overline 3$, $3 \overline 4$, $4 \overline 2$, $4 \overline 3$, and $4 \overline 4$.
See Figure~\ref{fig: D(G)} for an illustration.

Identifying $e_{\overline i}$ in $\RR^{\overline N}$ with $-e_{N+i}$ in $\RR^{2N}$ is a unimodular equivalence; thus, we have the following simple but important result.

\begin{lemma}\label{lem: ap is root}
	For all $G$, $\APQ_G$ is unimodularly equivalent to $Q_{D(G)}$. \qed
\end{lemma}

We now list two more theorems from \cite{Postnikov2009}.
In the first, $\sum$ denotes the Minkowski sum of polytopes and, given $S \subseteq [N]$, $\Delta_S = \conv\{ e_i \mid i \in S\}$.
Also, for a graph $G$ on $[N]$, set
\[
	P_{D(G)} = \sum_{i=1}^N \Delta_{\calN_G(i)\cup\{i\}} \subseteq  \RR^N.
\]
It is also written to reflect our particular context and does not quite capture the full strength of the original statement.
These two theorems are the last pieces needed to prove the main result of this section: Theorem~\ref{thm: translation}.

\begin{theorem}[{\cite[Theorem 12.2]{Postnikov2009}}]\label{thm: minkowski diff}
	Let $G$ be a graph on $[N]$ for which $D(G)$ is connected and let
	\[
		P^-_{D(G)} = \left\{x \in \RR^{N} \mid x + \Delta_{[N]} \subseteq P_{D(G)}\right\}.
	\]
	Then
	\[
		\NVol(Q_{D(G)}) = |P^-_{D(G)} \cap \ZZ^{N}|.
	\]
\end{theorem}

As written, Theorem~\ref{thm: minkowski diff} relies on $D(G)$ being connected.
Fortunately, the connectedness of $G$ is equivalent to the connectedness of $D(G)$.
We will use this fact occasionally so we present it as a lemma, although its proof is straightforward enough that we omit it.

\begin{lemma}\label{lem: connected iff}
	For any simple graph $G$, $G$ is connected if and only if $D(G)$ is connected. \qed
\end{lemma}

Since we are primarily working with $D(G)$ rather than $G$ directly, we let $\fkD(G)$ denote the set of $D(G)$-draconian sequences.

\begin{theorem}[{\cite[Theorem 11.3]{Postnikov2009}}]\label{thm: lattice pts draconian}
	Let $G$ be any graph.
	Then $|P^-_{D(G)} \cap \ZZ^{N}| = |\fkD(G)|$.
\end{theorem}

\begin{theorem}\label{thm: translation}
	For any connected graph $G$ on $[N]$, $\NVol(\APQ_G) = |\fkD(G)|$.
\end{theorem}

\begin{proof}
	Lemma~\ref{lem: connected iff} assures us that $D(G)$ is connected.
	By Lemma~\ref{lem: ap is root}, we know $\NVol(\APQ_G) = \NVol(Q_{D(G)})$.
	Applying Theorem~\ref{thm: minkowski diff} and Theorem~\ref{thm: lattice pts draconian} completes the proof. 
\end{proof}

To illustrate, let $G$ be the graph on $[4]$ with edges $12$, $23$, and $24$.
Here, we have $\calN_{D(G)}(1) = \{\overline 1,\overline 2\}$, $\calN_{D(G)}(2) = \{\overline 1,\overline 2,\overline 3,\overline 4\}$, $\calN_{D(G)}(3) = \{\overline 2,\overline 3\}$ and $\calN_{D(G)}(4) = \{\overline 2,\overline 4\}$.
Theorem~\ref{thm: translation} tells us that $\NVol(\APQ_G) = 8$ since 
\[\begin{aligned}
	\fkD(G) =& \{(0,3,0,0), (0,2,0,1), (1,1,1,0), (1,1,0,1), \\
		& \quad(1,0,1,1), (0,1,1,1), (0,2,1,0), (1,2,0,0)\}.
\end{aligned}\]

It will be very helpful for us to explicitly state when a sequence is $D(G)$-draconian.
The main difference is recognizing that for every vertex $i$ of $G$, $\deg_{D(G)}(i) = 1 + \deg_G(i)$. 

\begin{defn*}[Definition~\ref{def: draconian sequence}, rephrased]
	Let $G$ be a graph on $[N]$.
	Call $(a_1,\dots,a_N) \in \ZZ_{\geq 0}^N$ a \emph{$D(G)$-draconian sequence} if $\sum a_i = N-1$ and, for any $1 \leq i_1 < \cdots < i_k \leq N$,
	\[
		a_{i_1} + \cdots + a_{i_k} < \left|\bigcup_{j=1}^k \calN_{D(G)}(i_j)\right| = \left|\{i_1,\dots,i_k\} \cup \left(\bigcup_{j=1}^k \calN_G(i_j)\right)\right|
	\]
\end{defn*}

This translates our computation of normalized volume to a purely combinatorial computation.
The following simple observation will also be helpful at several points when proving the results in Section~\ref{sec: recurrences}.

\begin{remark}\label{rmk: relabeling invariant}
	The normalized volume of $\APQ_G$ is invariant under permutation of vertices.
\end{remark}

We now give a first nontrivial application of  Theorem~\ref{thm: translation} to an infinite class of graphs.

\begin{proposition}\label{prop: complete minus matching}
	Let $N>2$ and let $\mathcal{M}$ be any matching of size $k$ in $K_N$.
	Then
	\[
		\NVol(\APQ_{K_N \setminus \mathcal{M}}) = \binom{2(N-1)}{N-1} - 2k.
	\]
\end{proposition}

\begin{proof}
	Note that since $N>2$, $K_N \setminus \mathcal{M}$ is connected.
	First consider $k=0$.
	The $D(K_N)$-draconian sequences are the weak compositions of $N-1$ into $N$ parts, of which there are $\binom{2(N-1)}{N-1}$.
	When $k>0$, the deletion of each edge $uv$ in $M$ prohibits two compositions: those whose entries are all $0$ except for one, which is $N-1$ and located at position $u$ or $v$.
\end{proof}

Proposition~\ref{prop: complete minus matching} refers to a very specific class of graphs.
The next section proves results that allow for much more flexibility.


\section{Draconian recurrences}\label{sec: recurrences}

One of the main purposes of this article is to establish several recurrences for $\NVol(\APQ_G)$, using what we collectively call \emph{draconian recurrences}.
Certain specific recurrences will be given their own names as we encounter them.
For a simple first situation we consider the disjoint union of two graphs $G$ and $H$, which we denote $G+H$.
Since Theorem~\ref{thm: translation} only applies to connected graphs, we study their adjacency polytopes directly.

If $P \subseteq \RR^n$ and $Q\subseteq\RR^m$ are polytopes, each containing the origins $0_n$, $0_m$ respectively, then their \emph{free sum} is
\[
	P \oplus Q = \conv\{(P\times 0_m) \cup (0_n \times Q)\} \subseteq \RR^{n+m}.
\]
When $P$ and $Q$ are lattice polytopes, there is a convenient product formula we may invoke.

\begin{theorem}[{\cite[Theorem 2]{ChenDavis2019Product}}]\label{thm: free sum}
	Given full-dimensional convex polytopes $P \subseteq \RR^n$ and $Q \subseteq \RR^m$, if both $P$ and $Q$ contain the origin of their respective ambient spaces, then
	\[
		\NVol(P\oplus Q) = \NVol(P)\NVol(Q).
	\]
\end{theorem}

While Theorem~\ref{thm: free sum} insists that $P$ and $Q$ are full-dimensional, we may replace them with unimodularly equivalent polytopes $P' \subseteq \RR^{n'} \cong \aff(P)$ and $Q' \subseteq \RR^{m'} \cong \aff(P)$.
Since unimodular equivalence preserves normalized volume, the conclusion of Theorem~\ref{thm: free sum} remains true.
This gives us the last piece we need to prove the following.

\begin{proposition}\label{prop: disjoint union}
	If $G$ and $H$ are any two graphs, then
	\[
		\NVol(\APQ_{G+H}) = \NVol(\APQ_G)\NVol(\APQ_H).
	\]
\end{proposition}

\begin{proof}
	Let $|V(G)| = M$ and $|V(H)|=N$.
	First consider when $M=1$ or $N=1$.
	Without loss of generality we may assume $N=1$ and that the vertices of $G+H$ are labeled so that the isolated vertex from $H$ is labeled $M+1$.
	This means, in particular, that $\APQ_H$ consists of a single point, hence $\NVol(\APQ_H) = 1$.
	
	Let $A$ be the matrix whose columns are the vertices of $\APQ_G$.
	Partition $A$ as
	\[
		A = \left[\begin{array}{c}
			A_1 \\
			A_2
			\end{array}\right]
	\]
	where $A_1$ consists of the first $M$ rows of $A$ and $A_2$ consists of the last $M$ rows of $A$.
	The matrix of vertices of $\APQ_{G + H}$ can then be written as
	\[
		B = \left[\begin{array}{cc}
			A_1 & 0_{M \times 1}\\
			0_{1 \times \ell} & 1 \\
			A_2 & 0_{M \times 1} \\
			0_{1 \times \ell} & 1
			\end{array}\right]
	\]
	where $0_{k \times n}$ denotes the $k\times n$ all-zeros matrix and $\ell$ is the number of vertices of $\APQ_G$.
	It is clear that there is a unimodular transformation $f$ for which
	\[
		f(B) = \left[\begin{array}{cc}
			A_1 & 0_{M \times 1}\\
			0_{1 \times \ell} & 0 \\
			A_2 & 0_{M \times 1} \\
			0_{1 \times \ell} & 1
			\end{array}\right].
	\]
	Let $\pi(f(B))$ be the projection that drops row $M+1$ of $f(B)$, that is,
	\[
		\pi(f(B)) = \left[\begin{array}{cc}
			A_1 & 0_{M \times 1}\\
			A_2 & 0_{M \times 1} \\
			0_{1 \times \ell} & 1
			\end{array}\right].
	\]
	We now recognize $\pi(f(B))$ as a pyramid over $\APQ_G$.
	It is well-known that a lattice polytope and its pyramid have the same normalized volume (say, by implementing \cite[Theorem 2.4, Corollary 3.24]{BeckRobinsCCDed2}).
	Since $f$ is unimodular and since $\pi$ is a transformation providing a bijection between the lattices $\aff(f(B)) \cap \ZZ^{2(M+1)}$ and $\aff(\pi(f(B))) \cap \ZZ^{2M+1}$, we have
	\[
		\NVol(\APQ_{G+H}) = \NVol(\pi(f(B))) = \NVol(\APQ_G) = \NVol(\APQ_G)\NVol(\APQ_H).
	\]
	This proves the case of $M=1$ or $N=1$.

	Now assume $M,N \geq 2$.
	If $(x_1,\dots,x_{2M}) \in \APQ_G$, then, by construction,
	\[
	    \sum_{i=1}^M x_i = 1 \quad\text{and}\quad \sum_{i=M+1}^{2M} x_i = 1,
	\]
	and similar is true for $(y_1,\dots,y_{2N}) \in \APQ_H$.
	It follows that the polytopes
	\[
	    P = \{(x_2,\dots,x_M,x_{M+2},\dots,x_{2M}) \mid (x_1,\dots,x_{2M}) \in \APQ_G\}
	\]
	and
	\[
	    Q = \{(y_2,\dots,y_N,y_{N+2},\dots,y_{2N}) \mid (y_1,\dots,y_{2N}) \in \APQ_H\}
	\]
	are projections that are unimodularly equivalent to $\APQ_G$ and $\APQ_H$, respectively.
	Thus, $\NVol(\APQ_G) = \NVol(P)$ and $\NVol(\APQ_H) = \NVol(Q)$.
	Here, $P$ and $Q$ contain the origins of their respective ambient spaces, so 
	\[
		\NVol(P \oplus Q) = \NVol(P)\NVol(Q) = \NVol(\APQ_G)\NVol(\APQ_H).
	\]
	
	Label the vertices of $G+H$ using $[M+N]$ by adding $M$ to every vertex label of $H$.
	Let $f: \RR^{2M+2N} \to \RR^{2M+2N}$ be the map sending $(x_1,\dots,x_{2M+2N})$ to $(x_{\sigma(1)},\dots,x_{\sigma(2M+2N)})$ where
	\[
	    \sigma(i) = 
	        \begin{cases}
	            i & \text{ if } i \leq M \text{ or } i \geq 2M+N+1 \\
	            i+M & \text{ if } M+1 \leq i \leq M+N \\
	            i-N & \text{ if } M+N+1 \leq i \leq 2M+N.
	        \end{cases}
	\]
	Since $f$ only permutes coordinates it is a unimodular transformation.
	Moreover, the projection of $f(\APQ_{G+H})$ obtained from dropping the first, $(M+1)$th, $(2M+1)$th, and $(2M+N)$th coordinates is a lattice-preserving transformation sending $\APQ_{G+H}$ onto $P \oplus Q$.
	Therefore,
	\[
		\begin{aligned}
			\NVol(\APQ_{G+H}) &= \NVol(f(\APQ_{G+H})) \\
						&= \NVol(P \oplus Q) \\
						&= \NVol(P)\NVol(Q) \\
						&= \NVol(\APQ_G)\NVol(\APQ_H),
		\end{aligned}
	\]
	proving the result.
\end{proof}

In light of Proposition~\ref{prop: disjoint union}, we will focus for the rest of this section on graphs that are connected unless explicitly stated otherwise.
Restricting to when $G$ is connected allows us to use Theorem~\ref{thm: translation} and therefore we study the sets $\fkD(G)$ directly rather than relying on properties of their polytopes.

Recall that a graph $G$ is \emph{$k$-connected} if for any set $X$ of vertices, $|X| < k$, the subgraph $G - X$ is connected.
A \emph{cut-vertex} (respectively, \emph{cut-edge}) of $G$ is a vertex (respectively, edge) whose deletion from $G$ increases the number of components. 
A \emph{block} of a graph $G$ is an inclusion-maximal connected subgraph of $G$ with no cut-vertex.
Note that a block of a simple graph $G$ may be an isolated vertex, a cut-edge, or an inclusion-maximal $2$-connected subgraph of $G$.

\begin{theorem}\label{thm: cut vertices}
	Suppose $G$ is a connected graph with cut-vertex $v$ and $B$ is a block containing $v$.
	Setting $B' = G[\left(V(G) \setminus V(B)\right) \cup \{v\}]$ we have
	\[
		\NVol(\APQ_G) = \NVol(\APQ_B)\NVol(\APQ_{B'}).
	\]
\end{theorem}

\begin{proof}
	By Remark~\ref{rmk: relabeling invariant} we may assume without loss of generality that the cut-vertex is $1$, that $V(B)=[M]$, and that $V(B') = \{1,M+1,\dots,N\}$.
	We claim that the map
	\[
		f: \fkD(B) \times \fkD(B') \to \fkD(G)
	\]
	which sends $\left((c_1,c_2,\dots,c_M), (c'_1,c'_{M+1},\dots,c'_N) \right)$ to
	\[
		(d_1,\dots,d_N) = (c_1+c'_1,c_2,\dots,c_M,c'_{M+1},\dots,c'_N)
	\]
	is a well-defined bijection.

	For notational convenience set $c = (c_1,c_2,\dots,c_M)$ and $c' = (c'_1,c'_{M+1},\dots,c'_N)$.
	Since $c \in \fkD(B)$ and $c' \in \fkD(B')$, we know $\sum c_i = M-1$ and $\sum c'_i = N-M$.
	Thus, the sum of entries in $f(c,c')$ is $N-1$, one of the requirements for being $D(G)$-draconian.
	Now pick any sequence $1 \leq i_1 < \cdots < i_k \leq N$.
	If $i_k < M$ or $M < i_1$, then the corresponding $D(G)$-draconian inequality automatically holds. 
	So, suppose there is some positive $1 \leq \ell < k$ for which
	\[
		i_1 < \dots < i_\ell \leq M < i_{\ell+1} < \dots < i_k.
	\]
	If $1 < i_1$, then
	\[\begin{aligned}
		d_{i_1} + \cdots + d_{i_j} &= c_{i_1} + \cdots + c_{i_\ell} + c'_{i_\ell + 1} + \cdots + c'_k \\
			& < \left|\bigcup_{j=1}^{\ell} \calN_{D(B)}(i_j)\right| + \left|\bigcup_{j=\ell+1}^{k} \calN_{D(B')}(i_j)\right| - 1.
	\end{aligned}\]
	Since $B$ and $B'$ share just a single vertex, we have that
	\[
		\left|\bigcup_{j=1}^{\ell} \calN_{D(B)}(i_j)\right| + \left|\bigcup_{j=\ell+1}^k \calN_{D(B')}(i_j)\right| - 1 \leq \left|\bigcup_{j=1}^k \calN_{D(G)}(i_j)\right|.
	\]
	Chaining these inequalities together, the $D(G)$-draconian inequality holds.
	A similar argument holds if $1 = i_1$, only here we explicitly write $d_{i_1} = c_1 + c'_1$ and proceed as before.
	In both cases the $D(G)$-draconian inequality holds, therefore $f(c,c') \in \fkD(G)$.
	
	Showing that $f$ is injective is brief and straightforward, so we omit the details.
	What requires slightly more work is showing that $f$ is surjective.
	Let $d = (d_1,\dots,d_N) \in \fkD(G)$. 
	We claim that $d = f(c,c')$ where
	\[
	    c = \left(M-1-\sum_{i=2}^M d_i,d_2,\dots,d_M \right) \quad\text{and}\quad c' = \left(N-M-\sum_{j=M+1}^N d_j, d_{M+1},\dots, d_N\right)
	\]
	and $c \in \fkD(B)$, $c' \in \fkD(B')$.
	For notational convenience, we set
	\[
	    c_1 = M-1-\sum_{i=2}^M d_i \quad\text{and}\quad c'_1 = N-M-\sum_{j=M+1}^N d_j.
	\]
	Since it is clear that $d = f(c,c')$, the majority of the work will be in showing that $c \in \fkD(B)$ and $c' \in \fkD(B')$.
	The procedure is analogous for both, so we will only give the details for showing $c \in \fkD(B)$.
	
	By construction, the sum of entries in $c$ is $M-1$.
	Every inequality of the form
	\begin{equation}\label{eq: B}
	    d_{i_1} + \cdots + d_{i_k} < \left|\bigcup_{j=1}^k \calN_{D(B)}(i_j)\right|
	\end{equation}
	with $1 < i_1 < \cdots < i_k \leq M$ instantly holds since the neighbors of $2,\dots,M$ are the same in $D(G)$ and $D(B)$.
	It is also clear that $0 \leq c_1$ since, otherwise, $d_2 + \cdots + d_{M} > M-1$, which directly contradicts \eqref{eq: B}.
		
	Now consider a sum of a subsequence of $c$ of the form
	\[
	    c_1 + d_{i_1} + \cdots + d_{i_k}.
	\]
	By way of contradiction, suppose that this does not satisfy the corresponding $D(B)$-draconian inequality, that is,
	\[
	    c_1 + d_{i_1} + \cdots + d_{i_k} \geq \left|\calN_{D(B)}(1) \cup \left(\bigcup_{j=1}^k \calN_{D(B)}(i_j)\right)\right|.
	\]
	Since $1 < i_1 < i_k \leq M$, this inequality may be rewritten
	\begin{equation}\label{eq: surj1}
	    c_1 + d_{i_1} + \cdots + d_{i_k} \geq \left|\calN_{D(B)}(1) \cup \left(\bigcup_{j=1}^k \calN_{D(G)}(i_j)\right)\right|.
	\end{equation}
	We also now know that
	\begin{equation}\label{eq: surj2}
	    c'_1 + d_{M+1} + d_{M+2} + \cdots + d_{N} = N-M.
	\end{equation}
	Adding the corresponding sides of \eqref{eq: surj1} and \eqref{eq: surj2} and remembering that $c_1 + c'_1 = d_1$ results in 
	\[
	    d_1 + \sum_{j=1}^k d_{i_j} + \sum_{r = M+1}^{N} d_{r} \geq \left|\calN_{D(B)}(1) \cup \left(\bigcup_{j=1}^k \calN_{D(G)}(i_j)\right)\right| + N - M.
	\]
	Using the fact that $B'$ contains $N-M+1$ vertices,
	\[
	    d_1 + \sum_{j=1}^k d_{i_j} + \sum_{r = M+1}^{N} d_{r} \geq \left|\calN_{D(B)}(1) \cup \left(\bigcup_{j=1}^k \calN_{D(G)}(i_j)\right)\right| + \left|\bigcup_{r=M+1}^{N} \calN_{D(G)}(r)\right| - 1.
	\]
	Combining the first two summands on the right side counts the vertex $\overline{1}$ twice, resulting in
	\[
	    d_1 + \sum_{j=1}^k d_{i_j} + \sum_{r = M+1}^{N} d_{r} \geq \left|\calN_{D(B)}(1) \cup \left(\bigcup_{j=1}^k \calN_{D(G)}(i_j)\right) \cup \left(\bigcup_{r=M+1}^{N} \calN_{D(G)}(r)\right)\right|.
	\]
	which is a contradiction to $d$ being $D(G)$-draconian.
	Therefore the $D(B)$-draconian inequalities for $c$ all hold, and $c \in \fkD(B)$.
	An analogous argument shows $c' \in \fkD(B')$, proving $f$ is a bijection.
	This implies $|\fkD(B)||\fkD(B')| = |\fkD(G)|$; applying Theorem~\ref{thm: translation} completes the proof.
\end{proof}

The next result follows quickly from induction and the recurrences proven thus far.

\begin{corollary}\label{cor: forest volume}
	If $F$ is a forest on $N$ vertices with $k$ connected components, then we have $\NVol(\nabla_F^{\PQ}) = 2^{N-k}$. \qed
\end{corollary}

Interestingly, Corollary~\ref{cor: forest volume} implies that any two trees with the same number of edges will produce adjacency polytopes with the same normalized volume.
This does not happen for connected graphs in general: as we will show in Example~\ref{ex: cycles}, $\NVol(\APQ_{C_3}) = 6$, which is not the volume obtained from a path with three edges.
Moreover, even though two trees with the same number of vertices produce adjacency polytopes with equal normalized volumes, the polytopes themselves are not combinatorially equivalent.
Recall that the \emph{$f$-vector} of a polytope $P$ is the vector $(f_{-1},f_0,\dots,f_{\dim P})$ where $f_i$ is the number of $i$-dimensional faces of $P$, using the convention $f_{-1}=1$.

\begin{example}
    Let $G_1$ and $G_2$ be graphs on $[4]$.
    Let $E(G_1) = \{12, 23, 34\}$ and $E(G_2) = \{12, 13, 14\}$.
    One may verify the that the $f$-vector of $\APQ_{G_1}$ is 
    \[
        (1,10, 39, 77, 82, 46, 12, 1)
    \]
    and the $f$-vector of $\APQ_{G_2}$ is 
    \[
        (1,10, 39, 78, 86, 51, 14, 1).
    \]
    Thus the two polytopes are not combinatorially equivalent even though Theorem~\ref{thm: translation} guarantees that their normalized volumes are both $8$.
\end{example}
Through the recurrences established so far, we may reduce our work to considering only $2$-connected graphs.


\subsection{The subdivision recurrence}\label{subsec: subdivision}

Given $e \in E(G)$ let $G \col e$ denote the graph obtained by subdividing $e$.
Since we are using the convention $V(G) = [N]$, we will always assume that $V(G \col e) = [N+1]$.
The main result of this subsection is Theorem~\ref{thm: subdiv recurrence}, which gives a recurrence for $\NVol(\APQ_{G \col e})$ under certain conditions.
Establishing the recurrence requires multiple lemmas that have similar flavors but are distinct enough to warrant presenting their proofs.

The next three lemmas describe how to produce $D(G \col e)$-draconian sequences from $D(G)$-draconian sequences and $D(G\setminus e)$-draconian sequences.
We use the notation $A \uplus B$ to denote the disjoint union of the sets $A$ and $B$.

\begin{lemma}\label{lem: first rec lem}
	Let $G$ be any connected graph on $[N]$ and let $e=uv$ be an edge.
	If $c \in \fkD(G)$, then $\alpha(c) \in \fkD(G \col e)$ where $\alpha(c) = (c,1)$.
	Moreover, $\alpha$ is an injection.
\end{lemma}

\begin{proof}
	Let $c \in \fkD(G)$.
	By Remark~\ref{rmk: relabeling invariant} we may assume that $e = \{N-1,N\}$. 
	Showing that $\alpha$ is an injection is routine, so we focus mainly on showing $\alpha(c) \in \fkD(G \col e)$.
	
	Let $c = (c_1,\dots,c_{N})$.
	Since $c_1 + \cdots + c_{N} = N-1$, the sum of entries of $\alpha(c)$ is $N$.
	By construction, $\calN_{D(G)}(i) = \calN_{D(G \col e)}(i)$ for $i = 1,\dots,N-2$,
	\[
		\calN_{D(G \col e)}(N-1) = \left(\calN_{D(G)}(N-1)\setminus\{\overline{N}\}\right) \uplus \{\overline{N+1}\}
	\]
	and
	\[
		 \calN_{D(G \col e)}(N) = \left(\calN_{D(G)}(N)\setminus\{\overline{N-1}\}\right) \uplus \{\overline{N+1}\}.
	\]
	
	Pick a sequence $1 \leq i_1 < \cdots < i_k \leq N+1$.
	There are two cases to consider:
	\begin{enumerate}
		\item $\{\overline{N-1},\overline{N},\overline{N+1}\} \not\subseteq \bigcup_{j=1}^k \calN_{D(G \col e)}(i_j)$ and
		\item $\{\overline{N-1},\overline{N},\overline{N+1}\} \subseteq \bigcup_{j=1}^k \calN_{D(G \col e)}(i_j)$.
	\end{enumerate}
	In the first case, we can deduce two things: that $i_k \neq N+1$ and that if $N-1$ is one of the indices $i_1,\dots,i_k$, then no other neighbor of $N$ in $G$ is one of the indices $i_1,\dots,i_k$ (and vice versa).
	Therefore,
	\[
		\left|\bigcup_{j=1}^k \calN_{D(G)}(i_j)\right| = \left|\bigcup_{j=1}^k \calN_{D(G \col e)}(i_j)\right|
	\]
	and 
	\[
		c_{i_1} + \cdots + c_{i_k} < \left|\bigcup_{j=1}^k \calN_{D(G)}(i_j)\right| = \left|\bigcup_{j=1}^k \calN_{D(G \col e)}(i_j)\right|.
	\]
	In the second case, if $i_k < N+1$, we immediately get
	\[
		c_{i_1} + \cdots + c_{i_k}  < \left|\bigcup_{j=1}^k \calN_{D(G)}(i_j)\right| < \left|\bigcup_{j=1}^k \calN_{D(G \col e)}(i_j)\right|.
	\]
	Otherwise, $i_k = N+1$ and
	\[
		\left|\bigcup_{j=1}^k \calN_{D(G \col e)}(i_j)\right| \geq \left|\{\overline{N+1}\} \uplus \bigcup_{j=1}^{k-1} \calN_{D(G)}(i_j)\right| = \left|\bigcup_{j=1}^{k-1} \calN_{D(G)}(i_j)\right| + 1.
	\]
	This time, we get
	\[
		c_{i_1} + \cdots + c_{i_k}  = c_{i_1} + \cdots + c_{i_{k-1}}  + 1 < \left|\bigcup_{j=1}^{k-1} \calN_{D(G)}(i_j)\right| +1 \leq \left|\bigcup_{j=1}^k \calN_{D(G \col e)}(i_j)\right|.
	\]
	Since each case results in satisfying the $D(G \col e)$-draconian inequalities, we have shown that $\alpha(c) \in \fkD(G \col e)$.
\end{proof}

\begin{lemma}\label{lem: second rec lem}
	Let $G$ be a $2$-connected graph and let $e=uv$ be any edge.
	If $c \in \fkD(G\setminus e)$, then $\beta(c) \in \fkD(G \col e)$ where $\beta(c) = \alpha(c) + e_u - e_{N+1}$.
	Moreover, $\beta$ is an injection.
\end{lemma}
 
 \begin{proof}
    Arguing that $\beta$ is injective is routine, so its details are omitted.
    For what remains, by Remark~\ref{rmk: relabeling invariant} we may assume that $e = \{N-1,N\}$.
    We then want to show that, if $c = (c_1,\dots,c_{N}) \in \fkD(G \setminus e)$, then
    \[
        \beta(c) = (c_1,\dots,c_{N-2},c_{N-1}+1,c_{N},0) \in \fkD(G \col e).
    \]
    Note that, by Menger's theorem and the fact that $G$ is $2$-connected, $c$ exists since $G \setminus e$ is connected.
    
    Set $\beta(c) = (\beta_1,\dots,\beta_{N+1})$.
    Let $1 \leq i_1 < \cdots < i_k \leq N+1$ and set $\ell = k$ if $i_k < N+1$ and $\ell = k-1$ if $i_k = N+1$.
    If $N-1 \neq i_j$ for any $j$, then
    \[
        \beta_{i_1} + \cdots + \beta_{i_k} = c_{i_1} + \cdots + c_{i_\ell} < \left|\bigcup_{j=1}^\ell \calN_{D(G\setminus e)}(i_j)\right| \leq \left|\bigcup_{j=1}^k \calN_{D(G: e)}(i_j)\right|.
    \]
    Otherwise, $N-1 = i_j$ for some $j$. 
    In this case,
    \[
        \left|\bigcup_{j=1}^\ell \calN_{D(G\setminus e)}(i_j)\right| = \left|\bigcup_{j=1}^\ell \calN_{D(G \col e)}(i_j)\right| - 1 \leq \left|\bigcup_{j=1}^k \calN_{D(G \col e)}(i_j)\right| - 1.
    \]
    Together we have
    \[\begin{aligned}
        \beta_{i_1} + \cdots + \beta_{i_k} &= c_{i_1} + \cdots + c_{i_\ell} + 1 \\
        &< \left|\bigcup_{j=1}^\ell \calN_{D(G \setminus e)}(i_j)\right|+1\\
        &\leq \left|\bigcup_{j=1}^k \calN_{D(G \col e)}(i_j)\right|,
    \end{aligned}\]
    and the $D(G \col e)$-draconian inequality holds.
    Therefore $\beta(c) \in \fkD(G \col e)$.
\end{proof}

\begin{lemma}\label{lem: third rec lem}
	Let $G$ be a $2$-connected graph with an edge $e=uv$ such that $\deg_G(u)=2$ and the neighbors of $u$ are neighbors of each other.
	If $c \in \fkD(G)$, then $\gamma(c) \in \fkD(G \col e)$ where $\gamma(c)$ is formed by the following rule.
	Set $\gamma'(c) = \alpha(c) - e_u + e_{N+1}$.
	\begin{enumerate}
	    \item If $c \notin \fkD(G \setminus e)$, then
	    \begin{enumerate}
	        \item if $\gamma'(c) \in \fkD(G \col e)$, then set $\gamma(c) = \gamma'(c)$.
	        \item If $\gamma'(c) \notin \fkD(G \col e)$, then set $\gamma(c) = \alpha(c) + e_u - e_{N+1}$. 
	    \end{enumerate}
	    \item If $c \in \fkD(G \setminus e)$, then
	    \begin{enumerate}
	        \item if $\gamma'(c) \in \fkD(G \col e)$, then set $\gamma(c) = \gamma'(c)$.
	        \item If $\gamma'(c) \notin \fkD(G \col e)$, then set $\gamma(c) = \alpha(c) + e_v - e_{N+1}$.
	    \end{enumerate}
	\end{enumerate}
	Additionally, $\gamma$ is an injection.
\end{lemma}

\begin{proof}
    As usual, Remark~\ref{rmk: relabeling invariant} allows us to assume $e = \{N-1,N\}$ and $\deg_G(N-1) = 2$.
    This allows us to more specifically rewrite $\gamma$ as follows: set $\gamma'(c) = (c_1,\dots,c_{N-2},c_{N-1}-1,c_N,2)$.
    \begin{enumerate}
    \item If $c \notin \fkD(G \setminus e)$, then
        \begin{enumerate}
            \item if $\gamma'(c) \in \fkD(G \col e)$, then set $\gamma(c) = \gamma'(c)$.
            \item If $\gamma'(c) \notin \fkD(G \col e)$, then set $\gamma(c) = (c_1,\dots,c_{N-2},c_{N-1}+1,c_{N},0)$.
        \end{enumerate}
    \item If $c \in \fkD(G \setminus e)$, then
        \begin{enumerate}
            \item if $\gamma'(c) \in \fkD(G \col e)$, then set $\gamma(c) = \gamma'(c)$.
            \item If $\gamma'(c) \notin \fkD(G \col e)$, then set $\gamma(c) = (c_1,\dots,c_{N-1},c_{N}+1,0)$.
	    \end{enumerate}
	\end{enumerate}
	Throughout the proof we will use the notation
	\[
		\gamma'(c) = (\gamma'_1,\dots,\gamma'_{N+1}) \text{ and } \gamma(c) = (\gamma_1,\dots,\gamma_{N+1}).
	\]
    
    First suppose $c \notin \fkD(G \setminus e)$ and $\gamma'(c) \notin \fkD(G \col e)$, so that $\gamma(c) = (c_1,\dots,c_{N-2},c_{N-1}+1,c_{N},0)$.
    This places us in case $1$(b).
    Let $S$ be the set of all $c \in \fkD(G)$ satisfying these conditions.
    To show that $\gamma(c)$ is $D(G \col e)$-draconian we first show that $c_{N-1} \leq 1$.
    
    Partition $S$ into disjoint subsets $S_0 = \{ c \in S \mid c_{N-1}=0\}$ and $S_> = \{ c \in S \mid c_{N-1} > 0\}$.
    If $c \in S_0$, then clearly $c$ satisfies $c_{N-1} \leq 1$.
    For $c \in S_>$, since $c_{N-1} > 0$ we know $\gamma'_i \geq 0$ for all $i$, hence there must be some sequence $1 \leq i_1 < \dots < i_k \leq N+1$ for which $\gamma'(c)$ violates the corresponding $D(G \col e)$-draconian inequality.
    In fact, such a sequence cannot contain both $N-1$ and $N+1$, since, otherwise, the proof of Lemma~\ref{lem: first rec lem} implies
    \[
    	\gamma'_{i_1} + \cdots + \gamma'_{i_k} = c_{i_1} + \cdots + c_{i_{k-1}} + 1 < \left|\bigcup_{j=1}^k \calN_{D(G \col e)}(i_j)\right|
    \]
    and therefore the $D(G \col e)$-draconian inequality corresponding to $i_1 < \cdots < i_k$ is satisfied.
    Similarly, if $i_k < N+1$, then
    \[
    	\gamma'_{i_1} + \cdots + \gamma'_{i_k} \leq c_{i_1} + \cdots + c_{i_k},
    \]
    and the corresponding $D(G \col e)$-draconian inequality is satisfied again due to the proof of Lemma~\ref{lem: first rec lem}.
    Hence, any violation of a $D(G \col e)$-draconian inequality by $\gamma'(c)$ with $c \in S_>$ requires $i_k = N+1$ and $i_j \neq N-1$ for any $j < k$. 
    
    Since $\deg_G(N-1) = 2$, and since we are assuming $c_{N-1} > 0$, we only need to show that $c_{N-1} \neq 2$.
    If it were possible that $c_{N-1} = 2$, then for any sequence $1 \leq i_1 < \cdots < i_{k-1} < i_k = N+1$ not containing $N-1$, we would have
    \[
    	\begin{aligned}
		\gamma'_{i_1} + \cdots + \gamma'_{i_k} &= c_{i_1} + \cdots + c_{i_{k-1}} + 2 \\
									&= c_{i_1} + \cdots + c_{i_{k-1}} + c_{N-1} \\
									&< \left|\left(\bigcup_{j=1}^{k-1} \calN_{D(G)}(i_j)\right) \cup \calN_{D(G)}(N-1)\right|.
	\end{aligned}
    \]
    Notice that the set
    \[
    	T = \left(\bigcup_{j=1}^{k-1} \calN_{D(G)}(i_j)\right) \cup \calN_{D(G)}(N-1)
    \]
    contains both $\overline{N}$ and $\overline{N-1}$.
    On the other hand,
    \[
    	\bigcup_{j=1}^k \calN_{D(G \col e)}(i_j)
    \]
    contains $\overline{N+1}$ and all elements of $T$ with the potential exception of the vertex in $D(G \col e)$ corresponding to the neighbor of $N-1$ other than $N+1$, which is unique since we have assumed $\deg_G(N-1) = 2$.
    From this we conclude
    \[
    	\gamma'_{i_1} + \cdots + \gamma'_{i_k} < |T| \leq \left|\bigcup_{j=1}^k \calN_{D(G \col e)}(i_j)\right|.
    \]
    Since we have verified that all other $D(G \col e)$-draconian inequalities hold for $\gamma'(c)$, this would imply $\gamma'(c)$ is a $D(G \col e)$-draconian sequence, which is a contradiction.
    Therefore $c_{N-1} = 1$. 
    
    The remainder of our argument in establishing 1(b) applies to all elements of $S$.
    Consider a sum of the form $\gamma_{i_1} + \cdots + \gamma_{i_k}$ with $1 \leq i_1 < \cdots < i_k \leq N+1$.
    If $i_k < N-1$, then $\gamma_{i_j} = c_{i_j}$ and $\calN_{D(G)}(i_j) = \calN_{D(G \col e)}(i_j)$ for each $j = 1,\dots,k$, so that
    \[
    	\gamma_{i_1} + \cdots + \gamma_{i_k} = c_{i_1} + \cdots + c_{i_k} < \left|\bigcup_{j=1}^k \calN_{D(G)}(i_j)\right| = \left|\bigcup_{j=1}^k \calN_{D(G \col e)}(i_j)\right|.
    \]
    
    If $i_k = N-1$, then there are two subcases.
    First, if $\overline{N} \in \cup_{j=1}^{k-1} \calN_{D(G)}(i_j)$, then
    \[
    	\bigcup_{j=1}^k \calN_{D(G\col e)}(i_j) = \left(\bigcup_{j=1}^k \calN_{D(G)}(i_j)\right) \uplus \{\overline{N+1}\},
    \]
    hence
    \[
    	\gamma_{i_1} + \cdots + \gamma_{i_k} = c_{i_1} + \cdots + c_{i_k} + 1 < \left|\bigcup_{j=1}^k \calN_{D(G)}(i_j)\right| + \left| \{\overline{N+1}\}\right| = \left|\bigcup_{j=1}^k \calN_{D(G \col e)}(i_j)\right|.
    \]
    On the other hand, suppose $\overline{N} \notin \cup_{j=1}^{k-1} \calN_{D(G)}(i_j)$.
    Without loss of generality, we assume that the other neighbor of $N-1$ in $G$ is $N-2$.
    
    Since the neighbors of $N-1$ in $G$ are neighbors of each other, $i_{k-1} < N-2$, hence $\overline{N-1} \notin \cup_{j=1}^{k-1} \calN_{D(G)}(i_j)$ as well.
    Since we now know $c_{N-1} \leq 1$, we have
    \[
    	\begin{aligned}
		\gamma_{i_1} + \cdots + \gamma_{i_k} &= c_{i_1} + \cdots + c_{i_{k-1}} + c_{N-1} + 1 \\
			&\leq c_{i_1} + \cdots + c_{i_{k-1}} + 2 \\
			&< \left|\bigcup_{j=1}^{k-1} \calN_{D(G)}(i_j)\right| + \left| \{\overline{N-1}, \overline{N+1}\}\right| \\
			&\leq \left|\bigcup_{j=1}^k \calN_{D(G \col e)}(i_j)\right|.
	\end{aligned}
    \]
    This completes the case for $i_k = N-1$.

    If $i_k = N$ and $i_{k-1} < N-1$, then
    \[
    	\gamma_{i_1} + \cdots + \gamma_{i_k} = c_{i_1} + \cdots + c_{i_k} < \left|\bigcup_{j=1}^k \calN_{D(G)}(i_j)\right| \leq \left|\bigcup_{j=1}^k \calN_{D(G \col e)}(i_j)\right|,
    \]
    where the second inequality follows from recognizing that while $\overline{N-1}$ appears in the first union it may not appear in the second, and that $\overline{N+1}$ appears in the second union but does not appear in the first.
    If $i_k = N$ and $i_{k-1} = N-1$, then 
    \[
    	\gamma_{i_1} + \cdots + \gamma_{i_k} = c_{i_1} + \cdots + c_{i_k} + 1 < \left|\bigcup_{j=1}^k \calN_{D(G)}(i_j)\right| + \left|\{\overline{N+1}\}\right| = \left|\bigcup_{j=1}^k \calN_{D(G \col e)}(i_j)\right|.
    \]
    Finally, if $i_k = N+1$, then by the proof of Lemma~\ref{lem: first rec lem}, 
    \[
    	\gamma_{i_1} + \cdots + \gamma_{i_k} \leq c_{i_1} + \cdots + c_{i_{k-1}} + 1 < \left|\bigcup_{j=1}^k \calN_{D(G \col e)}(i_j)\right|.
    \]
    Therefore, $\gamma(c)$ is $D(G \col e)$-draconian when in case $1$(b).

    Now consider case $2$(b), so that $c \in \fkD(G \setminus e)$ and $\gamma'(c) \notin \fkD(G \col e)$.
    Proving that $\gamma(c) = (c_1,\dots,c_{N-1},c_{N}+1,0) \in \fkD(G \col e)$ follows the proof of Lemma~\ref{lem: second rec lem} almost identically, replacing $N-1$ with $N$.
    For this reason, we omit the details.

    To show that $\gamma$ is an injection, we can restrict to comparing the sequences of 1(a) with those of 2(a) and the sequences of 1(b) with those of 2(b).
    Fortunately, it is straightforward to see that no sequence can arise simultaneously as $\gamma(c)$ under the conditions of 1(a) and $\gamma(c')$ under the conditions of 2(a).
    If this were possible, we would obtain $c=c'$, but $c$ cannot simultaneously be a member of and absent from $\fkD(G\setminus e)$.
    
    For the remaining case, suppose $\gamma(c) = \gamma(d)$ where $\gamma(c)$ falls under the conditions of 1(b) and $\gamma(d)$ falls under the conditions of 2(b).
    Let $\gamma(d) = (\gamma''_1,\dots,\gamma''_{N+1})$.
    Since $\gamma(d)$ falls under the conditions of 2(b), we know $d \in \fkD(G \setminus e)$.
    In $G \setminus e$, the vertex $N-1$ has degree $1$, hence $\gamma''_{N-1} \in \{0,1\}$.
    We cannot have $\gamma''_{N-1} = 0$ since this would imply 
    \[
    	0 = \gamma''_{N-1} = \gamma_{N-1} = c_{N-1}+1,
    \]
    that is, we would have $c_{N-1} = -1$, which contradicts $c \in \fkD(G)$.
  
    If $\gamma''_{N-1} =  1$, then we claim $\gamma'(d) \in \fkD(G \col e)$ meaning $\gamma(d)$ would not be obtained from case 2(b).
    For ease of reference, we collect notation in terms of $\gamma''_1,\dots,\gamma''_{N+1}$ needed to complete the argument:
    \[
    	\begin{aligned}
		d &= (\gamma''_1,\dots,\gamma''_{N-2},1,\gamma''_{N}-1), \\
		c &= (\gamma''_1,\dots,\gamma''_{N-2},0,\gamma''_{N}), \\
		\gamma(d) &= (\gamma''_1,\dots,\gamma''_{N-2},1,\gamma''_{N},0), \\
		\gamma'(d) &= (\gamma''_1,\dots,\gamma''_{N-2},0,\gamma''_{N}-1,2). \\
	\end{aligned}
    \]
    
    Since $d \in \fkD(G \setminus e)$, all of the $D(G \col e)$-draconian inequalities involving the first $N$ coordinates of $\gamma'(d)$ are immediately satisfied.
    Consider, then, a sum involving the last coordinate of $\gamma'(d)$. 
    Note that $d \in \fkD(G \setminus e) \subseteq \fkD(G)$.
    If the sum involves the $N$th and the $(N-1)$th coordinates of $\gamma'(d)$ as well, then, since we know that $\alpha(d) \in \fkD(G \col e)$, we may write
    \[
    	\begin{aligned}
	    	\gamma''_{i_1} + \cdots + \gamma''_{i_{k-3}} + 0 + \gamma''_N-1 + 2 &= \gamma''_{i_1} + \cdots + \gamma''_{i_{k-3}} + 1 + (\gamma''_N-1) + 1 \\
		&< \left|\left(\bigcup_{j = 1}^{k-3} \calN_{D(G \col e)}(i_j)\right) \cup \left(\bigcup_{j=N-1}^{N+1} \calN_{D(G \col e)}(j)\right)\right|. \\
	\end{aligned}
    \]
    If the sum involves the $N$th coordinate of $\gamma'(d)$ but not the $(N-1)$th coordinate, then, this time, $\alpha(c) \in \fkD(G \col e)$ gives us
    \[
    	\begin{aligned}
	    	\gamma''_{i_1} + \cdots + \gamma''_{i_{k-2}} + \gamma''_N-1 + 2 &= \gamma''_{i_1} + \cdots + \gamma''_{i_{k-2}} + \gamma''_N + 1 \\
		&< \left|\left(\bigcup_{j = 1}^{k-2} \calN_{D(G \col e)}(i_j)\right) \cup \calN_{D(G \col e)}(N) \cup \calN_{D(G \col e)}(N+1)\right|. \\
	\end{aligned}
    \]
    A similar argument holds for when the sum involves the $(N-1)$th coordinate of $\gamma'(d)$ but not the $N$th coordinate.

    Next, suppose the sum is of the form 
    \[
    	\gamma''_{i_1} + \cdots + \gamma''_{i_{k-2}} + \gamma''_{N-2} + 2 = \gamma''_{i_1} + \cdots + \gamma''_{i_{k-2}} + \gamma''_{N-2} +1 + 1.
    \]
    Since $\alpha(d) \in \fkD(G \col e)$, we can say
    \[
    	\gamma''_{i_1} + \cdots + \gamma''_{i_{k-2}} + \gamma''_{N-2} +1 + 1 < \left|\left(\bigcup_{j = 1}^{k-2} \calN_{D(G \col e)}(i_j)\right) \cup \left(\bigcup_{j \in \{N-2,N-1,N+1\}} \calN_{D(G \col e)}(j)\right)\right|.
    \]
    Note that $\calN_{D(G \col e)}(N-1)$ can be dropped from this union because 
    \[
    	\calN_{D(G \col e)}(N-1) \subseteq \calN_{D(G \col e)}(N-2) \cup \calN_{D(G \col e)}(N+1).
    \]
    As a result, the desired inequality holds. 
    
    Finally, suppose the sum is of the form 
    \[
    	\gamma''_{i_1} + \cdots + \gamma''_{i_{k-1}} + 2
    \]
    where $i_{k-1} < N-2$.
    Since $d \in \fkD(G)$, we know
    \[
    	\gamma''_{i_1} + \cdots + \gamma''_{i_{k-1}} + 2 < \left|\bigcup_{j = 1}^{k-1} \calN_{D(G)}(i_j)\right| + 2.
    \]
    Since $i_{k-1} < N-2$, we know neither $\overline{N-1}$ nor $\overline{N+1}$, which are elements of $\calN_{D(G \col e)}(N+1)$, are in the above union. 
    Therefore,
    \[
    	\gamma''_{i_1} + \cdots + \gamma''_{i_{k-1}} + 2 < \left|\bigcup_{j = 1}^{k-1} \calN_{D(G)}(i_j)\right| + 2 \leq \left|\left(\bigcup_{j = 1}^{k-1} \calN_{D(G \col e)}(i_j)\right) \cup \calN_{D(G \col e)}(N+1)\right|.
    \]

    The above completes the argument that all $D(G \col e)$-draconian inequalities are satisfied by $\gamma'(d)$, i.e., $\gamma'(d) \in \fkD(G \col e)$.
    This contradicts the assumption that $\gamma(d)$ was obtained from case 2(b), so $\gamma''_{N-1} \neq 1$.
    Both possible values of $\gamma''_{N-1}$ lead to a contradiction, meaning no sequence obtained from case 1(b) can be obtained from case 2(b) and vice versa.
    Therefore, $\gamma$ is injective.
\end{proof}

Fix a particular edge $e$ of a $2$-connected graph $G$ for which one of the endpoints has degree $2$ in $G$.
Let $\scrA_G(e)$, $\scrB_G(e)$, and $\scrC_G(e)$ be the set of $D(G \col e)$-draconian sequences constructed from $\alpha$, $\beta$, and $\gamma$ in Lemmas~\ref{lem: first rec lem}, \ref{lem: second rec lem}, and \ref{lem: third rec lem}, respectively.

\begin{lemma}\label{lem: pairwise disjoint}
	Let $G$ be a $2$-connected graph with an edge $e = uv$ such that $\deg_G(u) = 2$ and the neighbors of $u$ are neighbors of each other.
	The sets $\scrA_G(e)$, $\scrB_G(e)$, and $\scrC_G(e)$ are pairwise disjoint.
\end{lemma}

\begin{proof}
    We continue to use the convention that $e = \{N-1,N\}$ and $\deg_G(N-1)=2$.
    We will also make the assumption that the other neighbor of $N-1$ in $G$ is $N-2$.
	By comparing the values of $c_{N+1}$, it is clear that $\scrA_G(e) \cap \scrB_G(e) = \emptyset$ and $\scrA_G(e) \cap \scrC_G(e) = \emptyset$.
	Thus we only need to focus on $\scrB_G(e) \cap \scrC_G(e)$.
	In fact, since $\gamma$ is an injection, we only need to consider elements of $\scrC_G(e)$ that fall under the conditions of 1(b) or 2(b) of Lemma~\ref{lem: third rec lem}.
	
	Suppose that $\beta(c) = \gamma(c')$ where
	\[
	    c = (c_1,\dots,c_N) \text{ and } c' = (c'_1, \dots, c'_N)
	\]
	are sequences with $c \in \fkD(G \setminus e)$ and $c' \in \fkD(G)$.
	If $\gamma(c')$ were to be constructed by the conditions of 1(b) in Lemma~\ref{lem: third rec lem}, it would follow that $c=c'$ since, in this case, $\beta(c) = \gamma(c')$ implies
	\[
		\alpha(c) + e_{N-1} - e_{N+1} = \alpha(c') + e_{N-1} - e_{N+1}
	\]
	and $\alpha$ is injective.
	However, this causes a contradiction, as Lemma~\ref{lem: second rec lem} requires $c \in \fkD(G \setminus e)$ while condition 1 of Lemma~\ref{lem: third rec lem} requires $c' = c \notin \fkD(G \setminus e)$.
	Hence we may assume $c \in \fkD(G \setminus e)$ and $\gamma(c')$ is constructed via condition 2(b) of Lemma~\ref{lem: third rec lem}.
	
	Since both $c,c' \in \fkD(G \setminus e)$, we know $c_{N-1},c'_{N-1} \leq 1$.
	By the definitions of $\beta$ and $\gamma$, we make several observations:
	\begin{itemize}
		\item $c_i = c'_i$ for all $i \leq N-2$;
		\item $\gamma(c')_{N-1} = c_{N-1} + 1$, hence $c'_{N-1} = 1$; and
		\item $\beta(c)_N = c'_N + 1$.
	\end{itemize}
	We claim that, in fact, $\gamma'(c') = (\gamma'_1,\dots,\gamma'_{N+1}) \in \fkD(G \col e)$, contradicting that $\gamma(c')$ was constructed via condition 2(b) of Lemma~\ref{lem: third rec lem}.
	Note that, based on our observations,
	\[
		\gamma'(c') = (c'_1,\dots,c'_{N-2}, c'_{N-1}-1, c'_N, 2) = (c_1, \dots ,c_{N-2}, 0, c_N-1,2).
	\]
	
	Consider a sum of the form $\gamma'_{i_1} + \cdots + \gamma'_{i_k}$ with $1 \leq i_1 < \cdots < i_k \leq N+1$.
	If $i_k \leq N-2$, then the neighbors of $i_j$ are the same in $D(G)$ and $D(G \col e)$, hence
	\[
		\gamma'_{i_1} + \cdots + \gamma'_{i_k} = c_{i_1} + \cdots + c_{i_k} < \left|\bigcup_{j=1}^k \calN_{D(G)}(i_j)\right| = \left|\bigcup_{j=1}^k \calN_{D(G \col e)}(i_j)\right|.
	\]
	If $i_k = N-1$, then 
	\[
		\gamma'_{i_1} + \cdots + \gamma'_{i_k} = c_{i_1} + \cdots + c_{i_{k-1}} + 0.
	\]
	As in the case when $i_k \leq N-2$, the neighbors of $i_1,\dots,i_{k-1}$ are the same in $D(G)$ and $D(G \col e)$, so that 
	\[
		\gamma'_{i_1} + \cdots + \gamma'_{i_k} < \left|\bigcup_{j=1}^{k-1} \calN_{D(G)}(i_j)\right| = \left|\bigcup_{j=1}^{k-1} \calN_{D(G:e)}(i_j)\right| \leq \left|\bigcup_{j=1}^k \calN_{D(G \col e)}(i_j)\right|.
	\]
	If $i_k = N$, then we have
	\[\begin{aligned}
		\gamma'_{i_1} + \cdots + \gamma'_{i_k} &= c_{i_1} + \cdots + c_{i_{k-1}} + c_N - 1 \\
			&< \left|\bigcup_{j=1}^k \calN_{D(G)}(i_j)\right| -1 \\
			&< \left|\bigcup_{j=1}^k \calN_{D(G \col e)}(i_j)\right|.
	\end{aligned}\]
	
	In the case $i_k = N+1$, we consider several subcases.
	If $i_{k-1} = N$, then
	\[
		\gamma'_{i_1} + \cdots + \gamma'_{i_k} = c_{i_1} + \cdots + c_{i_{k-1}} - 1 + 2 < \left|\bigcup_{j=1}^k \calN_{D(G \col e)}(i_j)\right|,
	\]
	where the inequality holds by Lemma~\ref{lem: first rec lem}.
	If $i_{k-1} = N-1$, then
	\[
		\gamma'_{i_1} + \cdots + \gamma'_{i_k} = c'_{i_1} + \cdots + c'_{i_{k-1}} - 1 + 2 < \left|\bigcup_{j=1}^k \calN_{D(G \col e)}(i_j)\right|,
	\]
	where the inequality again holds by Lemma~\ref{lem: first rec lem}.
	If $i_{k-1} = N-2$, then
	\[\begin{aligned}
		\gamma'_{i_1} + \cdots + \gamma'_{i_k} &= c'_{i_1} + \cdots + c'_{i_{k-1}} + c'_{N-1} - 1 + 2 \\
			&< \left|\left(\bigcup_{j=1}^k \calN_{D(G \col e)}(i_j)\right) \cup \calN_{D(G \col e)}(N-1) \right| \\
			&= \left|\bigcup_{j=1}^k \calN_{D(G \col e)}(i_j)\right|,
	\end{aligned}\]
	where the inequality follows from Lemma~\ref{lem: first rec lem} and the last equality comes from recognizing that $\calN_{D(G \col e)}(N-1) \subseteq \calN_{D(G \col e)}(N-2) \cup \calN_{D(G \col e)}(N+1)$, hence $\calN_{D(G \col e)}(N-1)$ may be freely dropped from the expression.
	Finally, if $i_{k-1} < N-2$, then 
	\[
		\gamma'_{i_1} + \cdots + \gamma'_{i_k} \leq c_{i_1} + \cdots + c_{i_{k-1}} + 2 < \left|\bigcup_{j=1}^k \calN_{D(G)}(i_j)\right| + \left|\{\overline{N-1},\overline{N+1}\}\right| \leq \left|\bigcup_{j=1}^k \calN_{D(G \col e)}(i_j)\right|.
	\]
	
	With this, we have verified that all of the $D(G \col e)$-draconian inequalities hold for $\gamma'(c')$, giving us the desired contradiction.
	Therefore, $\scrB_G(e) \cap \scrC_G(e) = \emptyset$, completing the proof.
\end{proof}

This result, together with the three lemmas preceding it, give $\scrA_G(e) \uplus \scrB_G(e) \uplus \scrC_G(e) \subseteq \fkD(G \col e)$.
It turns out that the reverse inclusion holds, establishing what we call \emph{the subdivision recurrence}.

\begin{theorem}[Subdivision recurrence]\label{thm: subdiv recurrence}
	Let $G$ be a $2$-connected graph with an edge $e=uv$ such that $\deg_G(u)=2$ and the neighbors of $u$ are neighbors of each other.
	Then $\fkD(G \col e) = \scrA_G(e) \uplus \scrB_G(e) \uplus \scrC_G(e)$ and, consequently,
	\[
		\NVol(\APQ_{G \col e}) = 2\NVol(\APQ_G) + \NVol(\APQ_{G \setminus e}).
	\]
\end{theorem}

\begin{proof}
    Again without loss of generality we may assume $e = \{N-1,N\}$ and $\deg_G(N-1)=2$.
	By Lemmas~\ref{lem: first rec lem}, \ref{lem: second rec lem}, and \ref{lem: third rec lem}, 
	\[
    	\scrA_G(e) \cup \scrB_G(e) \cup \scrC_G(e) \subseteq \fkD(G \col e).
    \]
	For the reverse inclusion, we will show that, given $d = (d_1,\dots,d_{N+1}) \in \fkD(G \col e)$, one of the following conditions holds:
	\begin{enumerate}
	    \item If $d_{N+1} = 2$, then $(d_1,\dots,d_{N-2},d_{N-1}+1,d_{N}) \in \fkD(G)$.
	    \item If $d_{N+1} = 1$, then $(d_1,\dots,d_{N}) \in \fkD(G)$.
	    \item If $d_{N+1} = 0$, then one of the following is true:
	        \begin{enumerate}
	            \item $(d_1,\dots,d_{N-2},d_{N-1}-1,d_{N}) \in \fkD(G \setminus e)$;
            	\item both $(d_1,\dots,d_{N-2},d_{N-1}-2,d_{N},2) \notin \fkD(G \col e)$ and $(d_1,\dots,d_{N-2},d_{N-1}-1,d_{N}) \in \fkD(G) \setminus \fkD(G \setminus e)$; or
	            \item both $(d_1,\dots,d_{N-2},d_{N-1}-1,d_{N}-1,2) \notin \fkD(G \col e)$ and $(d_1,\dots,d_{N-2},d_{N-1},d_{N}-1) \in \fkD(G \setminus e)$.
	        \end{enumerate}
	\end{enumerate}
	If the second condition holds, then $d \in \scrA_G(e)$; if condition 3(a) holds, then $d \in \scrB_G(e)$; if any of the remaining conditions hold, then $d \in \scrC_G(e)$.
	
	First suppose $d_{N+1}=2$ and let $1 \leq i_1 < \cdots < i_k \leq N$.
	Set $(c_1,\dots,c_{N}) = (d_1,\dots,d_{N-2},d_{N-1}+1,d_{N})$.
	If $i_k < N-1$, then $\calN_{D(G \col e)}(i_j) = \calN_{D(G)}(i_j)$ for each $j$, so the corresponding draconian inequality 
	\[
	    c_{i_1} + \cdots + c_{i_k} = d_{i_1} + \cdots + d_{i_k} < \left|\bigcup_{j=1}^k \calN_{D(G \col e)}(i_j)\right| = \left|\bigcup_{j=1}^k \calN_{D(G)}(i_j)\right|
	\]
	holds. 
	Otherwise, since $d \in \fkD(G \col e)$,
	\[\begin{aligned}
	    c_{i_1} + \cdots + c_{i_k} &\leq d_{i_1} + \cdots + d_{i_k} + 2 - 1 \\
	        &< \left|\left(\bigcup_{j=1}^k \calN_{D(G \col e)}(i_j)\right) \cup \calN_{D(G \col e)}(N+1)\right| - 1 \\
	        &= \left|\left(\bigcup_{j=1}^k \calN_{D(G)}(i_j)\right) \uplus \{\overline{N+1}\}\right| - 1 \\ 
	        &= \left|\bigcup_{j=1}^k \calN_{D(G)}(i_j)\right|.
	\end{aligned}\]
	Therefore, each $D(G)$-draconian inequality holds for $(c_1,\dots,c_N)$, establishing the first condition.
	
	Next suppose $d_{N+1} = 1$ and let $1 \leq i_1 < \cdots < i_k \leq N$.
	If $i_k < N-1$, then the corresponding draconian inequality holds as in the case of $d_{N+1}=2$. 
	If $i_k \geq N-1$, then we know from $d \in \fkD(G \col e)$ that
	\[\begin{aligned}
	    d_{i_1} + \cdots + d_{i_k} + 1 &< \left|\left(\bigcup_{j=1}^k \calN_{D(G \col e)}(i_j)\right) \cup \calN_{D(G \col e)}(N+1)\right| \\ 
	        &= \left|\left(\bigcup_{j=1}^k \calN_{D(G)}(i_j)\right) \uplus \{\overline{N+1}\}\right| \\ 
	        &= \left|\bigcup_{j=1}^k \calN_{D(G)}(i_j)\right| + 1.
	\end{aligned}\]
	Subtracting $1$ from both sides establishes the corresponding $D(G)$-draconian inequality for $(d_1,\dots,d_N)$.
	Thus the second condition holds.
	
	Establishing the last condition, where $d_{N+1} = 0$, requires the most care.
	Since $\deg_G(N-1) = 2$, we know that $d_{N-1} \in \{0,1,2\}$ and we will treat each case separately.
	
	Suppose $d_{N-1} = 0$.
	Our aim will be to show that condition 3(c) holds.
	It is clear that $(d_1,\dots,d_{N-2},d_{N-1}-1,d_{N}-1,2) \notin \fkD(G \col e)$ since $d_{N-1} - 1 < 0$.
	Now, if $d_{N}=0$, then there is a contradiction, since this and the $2$-connectivity of $G$ imply
	\[
	    N = d_1 + \cdots + d_{N-2} < \left|\bigcup_{j=1}^{N-2} \calN_{D(G \col e)}(j)\right| = \left|\bigcup_{j=1}^{N-2} \calN_{D(G)}(j)\right| = N.
	\]
	Thus, $d_{N} > 0$.
	
	Set $(c_1,\dots,c_{N}) = (d_1,\dots,d_{N-1},d_{N}-1)$ and consider the sum $c_{i_1} + \cdots + c_{i_k}$.
	If $i_k < N-1$, then the desired $D(G)$-draconian inequality holds using the same argument as for the previous conditions.
	If $i_k=N-1$, then
	\[
	    c_{i_1} + \cdots + c_{i_k} = d_{i_1} + \cdots + d_{i_{k-1}} 
	        < \left|\bigcup_{j=1}^{k-1} \calN_{D(G \col e)}(i_j)\right| 
	        = \left|\bigcup_{j=1}^{k-1} \calN_{D(G \setminus e)}(i_j)\right| 
	        \leq \left|\bigcup_{j=1}^k \calN_{D(G \setminus e)}(i_j)\right|.
	\]	

	Lastly, if $i_k = N$, then 
	\[\begin{aligned}
	    c_{i_1} + \cdots + c_{i_k} &= d_{i_1} + \cdots +d_{i_k}-1 \\
	        &< \left|\bigcup_{j=1}^k \calN_{D(G \col e)}(i_j)\right| - 1 \\
	        &= \left|\left(\bigcup_{j=1}^k \calN_{D(G \setminus e)}(i_j)\right) \uplus \{\overline{N+1}\}\right| - 1 \\ 
	        &= \left|\bigcup_{j=1}^k \calN_{D(G \setminus e)}(i_j)\right|.
	\end{aligned}\]
	Therefore, if $d_{N-1}=0$, then condition 3(c) holds.

	Next suppose $d_{N-1}=2$.
	Condition 3(c) clearly cannot hold since this condition requires $d_{N-1} \leq 1$, so we must show that either 3(a) or 3(b) holds.
	Suppose that condition 3(a) does not hold, that is, suppose $(d_1,\dots,d_{N-2},1,d_{N}) \notin \fkD(G \setminus e)$.
	Showing that this sequence is in $\fkD(G)$ can be done directly repeating our by-now-usual strategies, so the sequence is in $\fkD(G) \setminus \fkD(G \setminus e)$.

	To show that $(d_1,\dots,d_{N-2},0,d_{N},2) \notin \fkD(G \col e)$, observe that $(d_1,\dots,d_{N-2},1,d_{N}) \notin \fkD(G \setminus e)$ implies there is some inequality of the form
	\begin{equation}\label{eq: not 3a implies 3b}
	    d_{i_1} + \cdots + d_{i_k} \geq \left|\bigcup_{j=1}^k \calN_{D(G \setminus e)}(i_j)\right|
	\end{equation}
	with $i_k = N$ and $i_{k-1} < N-1$.
	If $\overline{N-1} \notin \bigcup_{j=1}^k \calN_{D(G \setminus e)}(i_j)$, then add $2$ to both sides of \eqref{eq: not 3a implies 3b} to get
	\[\begin{aligned}
	    d_{i_1} + \cdots + d_{i_k} + 2 &\geq \left|\bigcup_{j=1}^k \calN_{D(G \setminus e)}(i_j)\right| + 2 \\
	        &= \left|\bigcup_{j=1}^k \calN_{D(G \setminus e)}(i_j)\right| + \left|\{\overline{N-1},\overline{N+1}\}\right| \\
	        &= \left|\bigcup_{j=1}^k \calN_{D(G \col e)}(i_j)\right|,
	\end{aligned}\]
	which would imply $(d_1,\dots,d_{N-2},0,d_{N},2) \notin \fkD(G \col e)$.
	If $\overline{N-1} \in \bigcup_{j=1}^k \calN_{D(G \setminus e)}(i_j)$, then add $2$ to the left side of \eqref{eq: not 3a implies 3b} and $1 = \left|\{\overline{N}\}\right|$ to the right side; the conclusion is the same.
	Thus, if condition 3(a) does not hold, then condition 3(b) does hold. 
	
	For the case of when $d_{N-1}=1$, the first part of condition 3(b) clearly holds. 
	Verifying that $(d_1,\dots,d_{N-2},0,d_{N}) \in \fkD(G)$ is now routine, so either condition 3(a) holds or 3(b) holds.

	We have shown that, regardless of value of $d_{N+1}$, one of the three conditions holds, hence $d \in \scrA_G(e) \cup \scrB_G(e) \cup \scrC_G(e)$ and $\fkD(G \col e) = \scrA_G(e) \cup \scrB_G(e) \cup \scrC_G(e)$.
	By Lemma~\ref{lem: pairwise disjoint}, this union is disjoint, so
	\[\begin{aligned}
		|\fkD(G \col e)| &= |\scrA_G(e) \uplus \scrB_G(e) \uplus \scrC_G(e)| \\
			&= |\scrA_G(e)| + |\scrB_G(e)| + |\scrC_G(e)| \\
			&= 2|\fkD(G)| + |\fkD(G\setminus e)|.
	\end{aligned}\]
	Applying Theorem~\ref{thm: translation}, the result is proven.
\end{proof}

\begin{example}\label{ex: cycles}
	Consider $C_3 = ([3],\{12, 13, 23\})$ and let $e = 13$; there are six $D(C_3)$-draconian sequences:
	\[\begin{array}{ccc}
        		(2,0,0) & (0,2,0) & (0,0,2)  \\
		(1,1,0) & (1,0,1) & (0,1,1)
	\end{array}\]
	Subdividing $e$ replaces the edge $13$ with edges $34$ and $14$ to obtain $C_4$.
	By the subdivision recurrence, $\fkD(C_4) = \fkD(C_3 \col e) = \scrA_{C_3}(e) \uplus \scrB_{C_3}(e) \uplus \scrC_{C_3}(e)$.
	Following the definitions of $\alpha$, $\beta$, and $\gamma$ we obtain 
	\[\begin{aligned}
    		\scrA_{C_3}(e) &= \{(2,0,0,1), (0,2,0,1), (0,0,2,1), (1,1,0,1), (1,0,1,1), (0,1,1,1)\} \\
		\scrB_{C_3}(e) &= \{(1,2,0,0), (2,1,0,0), (2,0,1,0), (1,1,1,0)\} \\
		\scrC_{C_3}(e) &= \{(1,0,0,2),(1,0,2,0),(0,1,0,2),(0,0,1,2),(0,1,2,0),(0,2,1,0)\}. \\
	\end{aligned}\]
    Notice that $|\fkD(C_3)| = 3\cdot 2^1$ and $|\fkD(C_4)| = 4\cdot 2^2$.
\end{example}

Computational evidence suggests that the subdivision recurrence holds even without the condition that the neighbors of $u$ are neighbors of each other.
Indeed, Kohl~\cite{kohl} has verified that the recurrence holds for all graphs under this more relaxed condition having at most nine vertices.
However, the subdivision recurrence does not necessarily hold if we allow both endpoints of $e$ to have degree larger than $2$ in $G$.
For example, if $G = K_1 \vee P_3$, where $P_3$ is the path on three vertices, and $e$ is the edge of $G$ whose endpoints each have degree $3$ in $G$, then one may show that $\NVol(\APQ_{G \col e}) = 50$ whereas 
\[
	2\NVol(\APQ_{G}) + \NVol(\APQ_{G \setminus e}) = 2(18)+16 = 52.
\]

One important class of $2$-connected graphs that the subdivision recurrence does not directly cover is the class of cycles.
The conclusion of the subdivision recurrence holds, but to prove so requires a modified proof.

\begin{corollary}\label{cor: cycles}
	For a cycle $C_N$ on $N \geq 3$ vertices and any edge $e$ of $C_N$,
	\[
		\NVol(\APQ_{C_{N+1}}) = 2\NVol(\APQ_{C_N \col e}) + \NVol(\APQ_{C_N \setminus e}).
	\]
	Consequently, $\NVol(\APQ_{C_N}) = N2^{N-2}$.
\end{corollary}

\begin{proof}
	We will prove this result by showing that the conclusion of the subdivision recurrence holds for $C_N$ using the same functions $\alpha$, $\beta$, and $\gamma$ as in Lemmas~\ref{lem: first rec lem}, \ref{lem: second rec lem}, and \ref{lem: third rec lem}, respectively.
	Notice that none of the proofs of Lemmas~\ref{lem: first rec lem}, \ref{lem: second rec lem}, \ref{lem: pairwise disjoint}, or Theorem~\ref{thm: subdiv recurrence} rely on the neighbors of $u$ being neighbors of each other.
	In Lemma~\ref{lem: third rec lem}, the only time in which this condition is invoked is in establishing $\gamma(c) \in \fkD(G \col e)$ under case 1(b).
	Thus, by adapting the proof of case $1$(b) to $C_N$ we will have immediately established
	\[
		\NVol(\APQ_{C_{N+1}}) = 2\NVol(\APQ_{C_N}) + \NVol(\APQ_{C_N \setminus e}).
	\]
	The formula $\NVol(\APQ_{C_N}) = N2^{N-2}$ then follows from Corollary~\ref{cor: forest volume} and induction, whose details we omit.
	
	Without loss of generality we may assume that $C_N$ has vertex set $[N]$ and edges $\{1,2\}$, $\{2,3\}$,$\dots$,$\{N-1,N\}$,$\{1,N\}$, and we will subdivide the edge $e = \{N-1,N\}$.
	Let $c \in \fkD(C_N)$ so that $\gamma'(c) \notin \fkD(C_N \col e)$ and $c \notin \fkD(C_N \setminus e)$.
	Setting $c = (c_1,\dots,c_N)$, we have 
	\[
		\gamma(c) = (\gamma_1,\dots,\gamma_{N+1}) = (c_1,\dots,c_{N-2},c_{N-1}+1,c_N,0).
	\]

	Consider a sum of the form $\gamma_{i_1} + \cdots + \gamma_{i_k}$.
	If $i_k < N-1$, then $\calN_{D(C_N)}(i_j) = \calN_{D(C_N \col e)}(i_j)$ for each $j$, hence
	\[
		\gamma_{i_1} + \cdots + \gamma_{i_k} = c_{i_1} + \cdots + c_{i_k} < \left|\bigcup_{j=1}^k \calN_{D(C_N)}(i_j)\right| = \left|\bigcup_{j=1}^k \calN_{D(C_N \col e)}(i_j)\right|.
	\]
	
	When $i_{k} = N-1$ is when we have the most work to do.
	To begin, we may assume that $c_{i_j} > 0$ for all $j < k$ since, once the desired inequalities for these cases hold, if we were to check the inequality involving some $i_\ell < N-1$ for which $c_{i_\ell} = 0$, then we would instantly obtain
	\[
		\begin{aligned}
			\gamma_{i_1} + \cdots + \gamma_{i_k} + \gamma_{i_\ell} 
				&= \gamma_{i_1} + \cdots + \gamma_{i_k} \\
				&< \left|\bigcup_{j=1}^k \calN_{D(C_N \col e)}(i_j)\right| \\
				&\leq \left|\left(\bigcup_{j=1}^k \calN_{D(C_N \col e)}(i_j)\right) \cup \calN_{D(C_N \col e)}(i_\ell)\right|.
		\end{aligned}
	\]
	This allows us to assume for the rest of this case that $c_{i_j} > 0$ for all $j = 1,\dots,k-1$.
	
	If $i_1 = 1$, then
	\[
		\gamma_{i_1} + \cdots + \gamma_{i_k} = c_{i_1} + \cdots + c_{i_k} + 1 < \left|\bigcup_{j=1}^k \calN_{D(C_N)}(i_j)\right| + |\{\overline{N+1}\}| = \left|\bigcup_{j=1}^k \calN_{D(C_N \col e)}(i_j)\right|.
	\]
	If $i_1 > 1$, then partition the set $I = \{i_1,\dots,i_k\}$ into nonempty subsets $S_1,\dots,S_r$ such that each $S_i$ consists of consecutive integers and is maximal with that property with respect to containment.
	Additionally, label the subsets so that $\min(S_1) < \min(S_2) < \cdots < \min(S_r)$. 
	Define $t$ to be the largest index satisfying
	\[
		\min(S_{t}) - \max(S_{t-1}) > 2,
	\]
	or $t=1$ if no such index exists.
	Note that we always have $N-1 \in S_t \cup S_{t+1} \cup \cdots \cup S_r$.
	
	If $t>1$, then let $I_{<t} = S_1 \cup \cdots \cup S_{t-1}$ and $I_{\geq t} = S_t \cup \cdots \cup S_r$.
	Because $\min(I_{\geq t}) - \max(I_{<t}) > 2$, we know
	\[
		\left(\bigcup_{i_j \in I_{<t}} \calN_{D(G)}(i_j)\right) \cap \left(\bigcup_{i_j \in I_{\geq t}} \calN_{D(G)}(i_j)\right) = \emptyset.
	\]

	Since $c \in \fkD(C_N)$,
	\[
		\begin{aligned}
			c_{i_1} + \cdots + c_{i_k} &< \left|\bigcup_{i_j \in I_{<t}} \calN_{D(C_N)}(i_j)\right| + c_{\min(S_t)} + \cdots + c_{i_k} \\
				&< \left|\bigcup_{i_j \in I_{<t}} \calN_{D(C_N)}(i_j)\right| + \left|\bigcup_{i_j \in I_{\geq t}} \calN_{D(C_N)}(i_j)\right| \\
				&= \left|\bigcup_{j=1}^{k} \calN_{D(C_N)}(i_j)\right|.
		\end{aligned}
	\]
	Because of the two strict inequalities here, we may therefore say
	\[
		\begin{aligned}
			\gamma_{i_1} + \cdots + \gamma_{i_k} &= c_{i_1} + \cdots + c_{i_k} + 1 \\
				&\leq \left|\bigcup_{i_j \in I_{<t}} \calN_{D(C_N)}(i_j)\right| + c_{\min(S_t)} + \cdots + c_{i_k} \\
				&< \left|\bigcup_{i_j \in I_{<t}} \calN_{D(C_N)}(i_j)\right| + \left|\bigcup_{i_j \in I_{\geq t}} \calN_{D(C_N)}(i_j)\right| \\
				&= \left|\bigcup_{j=1}^{k} \calN_{D(C_N)}(i_j)\right|.
		\end{aligned}
	\]
		
	If $t = 1$, then there are two additional subcases to consider: when one of $c_1,c_N > 0$ and when $c_1 = c_N = 0$.
	In both subcases, we will use the fact that for each $\ell = 1, \dots, r$, 
	\[
		\sum_{i_j \in S_\ell} c_{i_j} \in \{|S_\ell|,|S_\ell| + 1\},
	\]
	that is, for a fixed $\ell$, the values of $c_{i_j}$ such that $i_j \in S_\ell$ are all $1$ or $2$ with at most one of them being $2$.
	Indeed, if there were two $2$s, then
	\[
		\sum_{i_j \in S_\ell} c_{i_j} \geq |S_\ell| + 2 = \left|\bigcup_{i_j \in S_\ell} \calN_{D(C_N)}(i_j)\right|,
	\]
	contradicting $c \in \fkD(C_N)$.
	We will moreover use fact that the case $t=1$ means
	\[
		\left|\bigcup_{j=1}^k \calN_{D(C_N)}(i_j)\right| = \left|\bigcup_{j=1}^k \calN_{D(C_N \col e)}(i_j)\right| = N - (i_1 - 2).
	\]
	
	Suppose that at least one of $c_1, c_N$ is positive.
	This implies there is at least one $\ell \in [r]$ for which $\sum_{i_j \in S_\ell} c_{i_j} = |S_\ell|$.
	Therefore,
	\[
		\begin{aligned}
			\gamma_{i_1} + \cdots + \gamma_{i_k} &= c_{i_1} + \cdots + c_{i_k} + 1 \\	
				&\leq \left(\sum_{i = 1}^r |S_i|+1\right) - 1 +1 \\
				&< N - (i_1 - 2) \\
				&= \left|\bigcup_{j=1}^k \calN_{D(C_N \col e)}(i_j)\right|.
		\end{aligned}
	\]
	
	If $c_1 = c_N = 0$, then we claim that $\gamma'(c') \in \fkD(C_N \col e)$, contradicting that we are in case 1(b).
	Here, $c_{N-1} > 0$ since, otherwise,
	\[
		c_2 + \cdots + c_{N-2} = N -1 = \left|\bigcup_{i = 2}^{N-2} \calN_{D(C_N)}(i)\right|,
	\]
	in which case $c \notin \fkD(C_N)$, a contradiction. 
	Thus, the entries of $\gamma'(c)$ are all nonnegative. 
	
	Let $\gamma'(c) = (\gamma'_1,\dots,\gamma'_{N+1})$ and consider a sum $\gamma'_{i_1} + \cdots + \gamma'_{i_k}$.
	If $i_k \leq N$, then
	\[
		\gamma'_{i_1} + \cdots + \gamma'_{i_k} \leq c_{i_1} + \cdots + c_{i_k} < \left|\bigcup_{j=1}^k \calN_{D(C_N)}(i_j)\right| \leq \left|\bigcup_{j=1}^k \calN_{D(C_N \col e)}(i_j)\right|.
	\]
	Otherwise, $i_k = N+1$.
	Since we are in the case $c_1 = c_N = 0$, we may assume $i_1 > 1$ and $i_{k-1} < N$.
	if $i_{k-1} = N-1$, then $\alpha(c) \in \fkD(C_N \col e)$ gives us
	\[
		\gamma'_{i_1} + \cdots + \gamma'_{i_k} = c_{i_1} + \cdots + c_{i_{k-1}} - 1 + 2 < \left|\bigcup_{j=1}^k \calN_{D(C_N \col e)}(i_j)\right|
	\]
	right away.
	If $i_{k-1} < N-1$, then neither $\overline{N}$ nor $\overline{N+1}$ is a neighbor of any of $i_1,\dots,i_{k-1}$ in either $D(C_N)$ or $D(C_N \col e)$, so
	\[
		\gamma'_{i_1} + \cdots + \gamma'_{i_k} = c_{i_1} + \cdots + c_{i_{k-1}} + 2 < \left|\bigcup_{j=1}^k \calN_{D(C_N)}(i_j)\right| +  |\{\overline{N}, \overline{N+1}\}| \leq \left|\bigcup_{j=1}^k \calN_{D(C_N \col e)}(i_j)\right|.
	\]
	Therefore, $\gamma'(c) \in \fkD(C_N \col e)$, contradicting the fact that $\gamma(c)$ arises from case 1(b).
	This completes the case $i_k = N-1$.
	
	If $i_k = N$, then consider $i_{k-1}$.
	If $i_{k-1} < N-1$, then we obtain the desired inequality again from knowing $\alpha(c) \in \fkD(C_N \col e)$.
	Otherwise $i_{k-1} = N-1$, in which case
	\[
		\begin{aligned}
			\gamma_{i_1} + \cdots + \gamma_{i_k} &= c_{i_1} + \cdots + c_{i_k} + 1\\
				&< \left|\bigcup_{j=1}^k \calN_{D(C_N)}(i_j)\right| + |\{\overline{N+1}\}| \\
				&= \left|\bigcup_{j=1}^k \calN_{D(C_N \col e)}(i_j)\right|.
		\end{aligned}
	\]

	Lastly, if $i_k = N+1$, then $\gamma_{i_k} = 0$, and by the previous cases we have
	\[
		\gamma_{i_1} + \cdots + \gamma_{i_k} = \gamma_{i_1} + \cdots + \gamma_{i_{k-1}} < \left|\bigcup_{j=1}^{k-1} \calN_{D(C_N \col e)}(i_j)\right| \leq \left|\bigcup_{j=1}^k \calN_{D(C_N \col e)}(i_j)\right|.
	\]
	This final inequality establishes $\gamma(c) \in \fkD(C_N \col e)$ for the case 1(b), which completes our proof. 
\end{proof}

We close this subsection with an invitation to the reader.

\begin{question}
    Under what conditions for a graph $G$ and an edge $e$ is there a ``nice'' recurrence for $\NVol(\APQ_{G \col e})$?
\end{question}


\subsection{The triangle recurrence}\label{subsec: triangle}

The framework which establishes the subdivision recurrence can be adapted to a different operation.
Given an edge $e = uv$ of a graph $G$, let $G \triangle e$ denote the graph with edge set $E(G) \cup \{uw,vw\}$ where $w$ is a new vertex.
We will continue to assume $V(G) = [N]$ and $V(G\triangle e) = [N+1]$.
As in Section~\ref{subsec: subdivision}, establishing a recurrence formula for $\fkD(G\triangle e)$ will require establishing several smaller results first.
The first two of these have proofs analogous enough to the proofs of Lemma~\ref{lem: first rec lem} and Lemma~\ref{lem: second rec lem}, respectively, that we omit their details.

\begin{lemma}\label{lem: first rec lem tri}
    Let $G$ be any connected graph on $[N]$ and $e$ any edge.
	If $c \in \fkD(G)$, then $\alphatri(c) \in \fkD(G\triangle e)$ where $\alphatri(c) = (c,1)$.
	Moreover, $\alphatri$ is injective. \qed
\end{lemma}

\begin{lemma}\label{lem: second rec lem tri}
        Let $G$ be a connected graph on $[N]$ and let $e = uv$ be any edge.
	If $c \in \fkD(G)$, then $\betatri(c) \in \fkD(G\triangle e)$ where
	\[
		\betatri(c) = \alphatri(c) + e_u - e_{N+1}.
	\]
	Additionally, $\betatri$ is injective. \qed
\end{lemma}

The next lemma is analogous to Lemmas~\ref{lem: third rec lem}, but this time its proof is different enough for us to justify providing it.
It will be helpful to introduce the analogues of $\scrA_G(e)$ and $\scrB_G(e)$ here: let $\scrAtri_G(e)$ and $\scrBtri_G(e)$ be the $D(G\triangle e)$-draconian sequences constructed with $\alphatri$ and $\betatri$ in Lemmas~\ref{lem: first rec lem tri} and \ref{lem: second rec lem tri}, respectively.

\begin{lemma}\label{lem: third rec lem tri}
    Let $G$ be a connected graph on $[N]$ and let $e = uv$ be any edge for which $\deg_G(u) = 2$.
	If $c \in \fkD(G)$, then $\gammatri(c) \in \fkD(G\triangle e)$ where
	\[
		\gammatri(c) =
		    \begin{cases}
		        \alphatri(c) + e_v - e_{N+1} & \text{ if not in } \scrBtri_G(e) \\
		        \alphatri(c) - e_u + e_{N+1} & \text{ otherwise}.
		    \end{cases}
	\]
	Additionally, $\gammatri$ is injective.
\end{lemma}

\begin{proof}
    That $\gammatri$ is injective is clear.
	For what remains, by Remark~\ref{rmk: relabeling invariant} we again assume without loss of generality that $e = \{N-1,N\}$.
	We also assume that the other neighbor of $N-1$ in $G$ is $N-2$.
	So, if $c = (c_1,\dots,c_{N}) \in \fkD(G)$, then we must prove $\gammatri(c) \in \fkD(G\triangle e)$ where
	\[
		\gammatri(c) = (\gammatri_1,\dots,\gammatri_{N+1}) = 
		\begin{cases}
		(c_1,\dots,c_{N-2},c_{N-1},c_{N}+1,0) & \text{ if not in } \scrBtri_G(e) \\
		(c_1,\dots,c_{N-2},c_{N-1}-1,c_{N},2) & \text{ otherwise}.
		\end{cases}
	\]
	Note that, in both cases, the entries sum to $N$.
	
	If $\gammatri(c) = (c_1,\dots,c_{N-2},c_{N-1},c_{N}+1,0)$, then showing it is $D(G \triangle e)$-draconian is entirely analogous to the proof of Lemma~\ref{lem: second rec lem}. 
	Otherwise, $\gammatri(c) = (c_1,\dots,c_{N-2},c_{N-1}-1,c_{N},2)$.
	Being in this case means that $(c_1,\dots,c_{N-1},c_N+1,0) \in \scrBtri(e)$.
	Hence, we know $c_{N-1} \geq 1$, so all entries of $\gammatri$ are nonnegative, and we also know
	\begin{equation}\label{eq: c'}
		c' = (c_1,\dots, c_{N-2}, c_{N-1}-1,c_N+1) \in \fkD(G).
	\end{equation}
	
	Consider a sum $\gammatri_{i_1} + \cdots + \gammatri_{i_k}$ with $1 \leq i_1 < \dots < i_k \leq N+1$.
	If $i_k \leq N$, then
	\[
		\gammatri_{i_1} + \cdots + \gammatri_{i_k} \leq c_{i_1} + \cdots + c_{i_k} < \left|\bigcup_{j=1}^k \calN_{D(G)}(i_j)\right| \leq \left|\bigcup_{j=1}^k \calN_{D(G \triangle e)}(i_j)\right|.
	\]
	If $i_k = N+1$ then there are four subcases to consider depending on the value of $i_{k-1}$.
	
	If $i_{k-1} = N$, then we may write the sum as
	\[
		\gammatri_{i_1} + \cdots + \gammatri_{i_k} = c_{i_1} + \cdots + c_{i_{k-2}} + (c_N + 1) + 1.
	\]
	By \eqref{eq: c'}, we know that $\alphatri(c') \in \fkD(G \triangle e)$.
	The above sum appears when verifying this fact, so we know that
	\[
		\gammatri_{i_1} + \cdots + \gammatri_{i_k} < \left|\bigcup_{j=1}^k \calN_{D(G \triangle e)}(i_j)\right|.
	\]
	If $i_{k-1} = N-1$, then 
	\[
		\gammatri_{i_1} + \cdots + \gammatri_{i_k} = c_{i_1} + \cdots + c_{i_{k-2}} + c_{N-1} + 1 < \left|\bigcup_{j=1}^k \calN_{D(G \triangle e)}(i_j)\right|,
	\]
	where the inequality again comes from knowing $\alphatri(c) \in \fkD(G \triangle e)$.
	If $i_{k-1} = N-2$, then again by applying $\alphatri$, we may say
	\[
		\gammatri_{i_1} + \cdots + \gammatri_{i_k} \leq c_{i_1} + \cdots + c_{i_{k-1}} + c_{N-1} - 1 + 2  < \left|\left(\bigcup_{j=1}^k \calN_{D(G \triangle e)}(i_j)\right) \cup \calN_{D(G \triangle e)}(N-1) \right|.
	\]
	Noticing that 
	\[
		\calN_{D(G \triangle e)}(N-1) \subseteq \calN_{D(G \triangle e)}(N-2) \cup \calN_{D(G \triangle e)}(N+1) \subseteq \bigcup_{j=1}^k \calN_{D(G \triangle e)}(i_j)
	\]
	we may drop $\calN_{D(G \triangle e)}(N-1)$ from the union, which establishes the desired inequality.
	Lastly, if $i_{k-1} < N-2$, then neither $\overline{N-1}$ nor $\overline{N+1}$ is a neighbor of $i_j$ in $D(G \triangle e)$ for $j \leq k-1$, so we may say
	\[
		\gammatri_{i_1} + \cdots + \gammatri_{i_k} < \left|\bigcup_{j=1}^{k-1} \calN_{D(G \triangle e)}(i_j) \uplus \{\overline{N-1},\overline{N+1}\}\right| \leq \left|\bigcup_{j=1}^k \calN_{D(G \triangle e)}(i_j)\right|.
	\]
	In all cases, the required $D(G \triangle e)$-draconian inequality holds.
	Therefore, we have shown $\gammatri(c) \in \fkD(G \triangle e)$ for all $c \in \fkD(G)$.
\end{proof}

Let $\scrCtri_G(e)$ be the $D(G\triangle e)$-draconian sequences constructed from $\gammatri$ in Lemma~\ref{lem: third rec lem tri}.
The proof of the following is completely analogous to the proof of Lemma~\ref{lem: pairwise disjoint}.

\begin{lemma}\label{lem: pairwise disjoint tri}
    Let $G$ be a graph having an edge $e=uv$ with $\deg_G(u)=2$.
    The sets $\scrAtri_G(e)$, $\scrBtri_G(e)$, and $\scrCtri_G(e)$ are pairwise disjoint. \qed
\end{lemma}

As in Section~\ref{subsec: subdivision}, the previous four lemmas imply $\scrAtri_G(e) \uplus \scrBtri_G(e) \uplus \scrCtri_G(e) \subseteq \fkD(G \triangle e)$.
The reverse inclusion again holds under certain restrictions, establishing what we call the \emph{triangle recurrence}..
We present the proof below, deferring portions of it to two lemmas afterward.

\begin{theorem}[Triangle Recurrence]\label{thm: triangle recurrence}
	Let $G$ be any connected graph and let $e = uv$ be an edge with $\deg_{G}(u) = 2$.
	If $\deg_G(v) = 2$ or if the neighbors of $u$ are neighbors of each other, then
	\[
		\NVol(\APQ_{G \triangle e}) = 3\NVol(\APQ_G).
	\]
\end{theorem}

\begin{proof}
    As usual we assume $V(G) = [N]$, $e = \{N-1,N\}$, and $\deg_{G}(N-1) = 2$.
    We will further assume that the other neighbor of $N-1$ in $G$ is $N-2$.
    Lemmas~\ref{lem: first rec lem tri}, \ref{lem: second rec lem tri}, and \ref{lem: third rec lem tri} show that 
    \[
        \scrAtri_G(e) \cup \scrBtri_G(e) \cup \scrCtri_G(e) \subseteq \fkD(G\triangle e),
    \]
    so we must show the reverse inclusion holds.
    
    Let $d = (d_1,\dots,d_{N+1}) \in \fkD(G \triangle e)$.
    As with the subdivision recurrence, there are three statements we must establish:
    \begin{itemize}
        \item If $d_{N+1} = 0$, then $(d_1,\dots,d_{N-2},d_{N-1}-1,d_{N}) \in \fkD(G)$ or, if this is not the case, then $(d_1,\dots,d_{N-2},d_{N-1},d_{N}-1) \in \fkD(G)$;
       	\item If $d_{N+1} = 1$, then $(d_1,\dots,d_{N}) \in \fkD(G)$; and
        \item If $d_{N+1} = 2$, then both $(d_1,\dots,d_{N-2},d_{N-1}+1,d_{N}+1,0) \in \scrBtri_G(e)$ as well as $(d_1,\dots,d_{N-2},d_{N-1}+1,d_{N}) \in \fkD(G)$.
    \end{itemize}
    For ease of readability, the case $d_{N+1} = 0$ is deferred to Lemmas~\ref{lem: partial reverse 1} and \ref{lem: partial reverse 2}, where the two different conditions on the vertices $N-1$ and $N$ are treated individually.
    
    Suppose, then, that $d_{N+1} = 1$.
    Pick any $D(G \triangle e)$-draconian sequence of the form $(d_1,\dots,d_{N},1)$.
    Let $1 \leq i_1 < \cdots < i_k \leq N$. 
    If $i_j \neq N-1, N$ for all $j$, then the neighbors of $i_j$ are the same in $D(G \triangle e)$ and $D(G)$, so the corresponding $D(G)$-draconian inequality instantly holds.
    Otherwise,
    \[\begin{aligned}
        d_{i_1} + \cdots + d_{i_k} &= d_{i_1} + \cdots + d_{i_k} + 1 - 1\\
        &< \left|\left(\bigcup_{j=1}^k \calN_{D(G \triangle e)}(i_j)\right) \cup \calN_{D(G \triangle e)}(N+1)\right| - 1 \\
        &= \left|\left(\bigcup_{j=1}^k \calN_{D(G \triangle e)}(i_j)\right) \cup \left\{\overline{N+1}\right\} \right| - 1 \\
        &= \left|\bigcup_{j=1}^k \calN_{D(G)}(i_j)\right|
    \end{aligned}\]
    Thus $(d_1,\dots,d_{N}) \in \fkD(G)$.
    
    Lastly, suppose $d_{N+1} = 2$. 
    For this case we first show that $(d_1,\dots,d_{N-1}+1, d_{N}+1,0) \in \scrBtri_G(e)$.
    This can be rephrased as wanting to show $(d_1,\dots,d_{N-1}+1, d_{N}+1,0) = \betatri(c)$ for some $c$, or, in yet other words, that $(d_1,\dots,d_{N-1}, d_{N}+1) \in \fkD(G)$.
    
    Set $d' = (d'_1,\dots,d'_{N}) = (d_1,\dots,d_{N-1}, d_{N}+1)$ and consider $1 \leq i_1 < \cdots < i_k \leq N$. 
    If $i_k < N-1$, then the $D(G)$-draconian inequality holds as usual. 
    If $i_k = N-1$, then observe
    \[\begin{aligned}
        d'_{i_1} + \cdots + d'_{i_k}
            &< d_{i_1} + \cdots + d_{i_k} + 2 - 1 \\ 
            &<  \left|\left(\bigcup_{j=1}^k \calN_{D(G)}(i_j)\right) \cup \{\overline{N+1}\}\right| - 1 \\
            &= \left|\bigcup_{j=1}^k \calN_{D(G)}(i_j)\right|.
    \end{aligned}\]
    If $i_k = N$, then repeat this argument but by inserting ``$+1-1$'' instead of ``$+2-1$''.
    In all cases, the $D(G)$-draconian inequality holds, so $d' \in \fkD(G)$, as needed.
    In fact, showing that $(d_1,\dots,d_{N-2},d_{N-1}+1,d_{N})$ is $D(G)$-draconian has an entirely analogous argument. 
    Therefore, this completes the case for $d_{N+1}=2$.
    
	By Lemma~\ref{lem: pairwise disjoint tri},
	\[  
	    \fkD(G\triangle e) = \scrAtri_G(e) \uplus \scrBtri_G(e) \uplus \scrCtri_G(e).
	\]
	Thus,
	\[
	   |\fkD(G\triangle e)| = |\scrAtri_G(e)| + |\scrBtri_G(e)| + |\scrCtri_G(e)| = 3|\fkD(G)|.
	\]
	Finally, by Theorem~\ref{thm: translation}, we obtain
	\[
		\NVol(\APQ_{G \triangle e}) = 3\NVol(\APQ_G).
	\]
\end{proof}

\begin{example}
	Let $C_3$ be the $3$-cycle as in Example~\ref{ex: cycles} and again choose $e = 13$.
	The $D(G \triangle e)$-draconian sequences are formed from the disjoint union of the three sets
	\[\begin{aligned}
		\scrAtri_{C_3}(e) &= \{(2,0,0,1),(0,2,0,1),(0,0,2,1),(1,1,0,1),(1,0,1,1),(0,1,1,1)\}\\
		\scrBtri_{C_3}(e) &= \{(3,0,0,0),(1,2,0,0),(1,0,2,0),(2,1,0,0),(2,0,1,0),(1,1,1,0)\} \\
		\scrCtri_{C_3}(e) &= \{(1,0,0,2),(0,2,1,0),(0,0,3,0),(0,1,0,2),(0,0,1,2),(0,1,2,0)\} \\
	\end{aligned}\]
\end{example}

As in the case of the subdivision recurrence, by relaxing the requirement that $e$ has an endpoint of degree $2$ in $G$, the result may no longer hold.
The same example as before, where $G = K_1 \vee P_3$ and $e$ is the edge whose endpoints each have degree $3$ in $G$, demonstrates this.
The normalized volume of $\APQ_{G\triangle e}$ is $52$ whereas a naive attempt to apply the triangle recurrence would predict $54$.

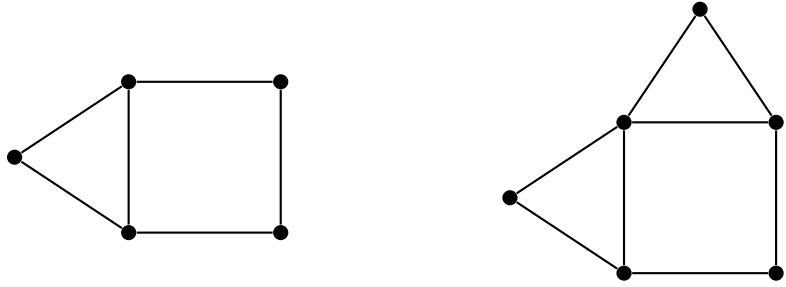
\begin{figure}
\begin{center}
\begin{tikzpicture}
\begin{scope}[every node/.style={circle,fill,inner sep=0pt,minimum size=2mm}]
	\node (A) at (0,0) {};
	\node (B) at (0,2) {};
	\node (C) at (2,2) {};
	\node (D) at (2,0) {};
	\node (E) at (-1.5,1) {};
\end{scope}
\node (F) at (0,-0.5) {};

\draw[thick] (A) -- (B) -- (C) -- (D) -- (A) -- (E) -- (B);
\end{tikzpicture}
\hspace{1in}
\begin{tikzpicture}
\begin{scope}[every node/.style={circle,fill,inner sep=0pt,minimum size=2mm}]
	\node (A) at (0,0) {};
	\node (B) at (0,2) {};
	\node (C) at (2,2) {};
	\node (D) at (2,0) {};
	\node (E) at (-1.5,1) {};
	\node (F) at (1,3.5) {};
\end{scope}

\draw[thick] (A) -- (B) -- (C) -- (D) -- (A) -- (E) -- (B) -- (F) -- (C);
\end{tikzpicture}
\end{center}
\caption{Two graphs $G_1$ (left) and $G_2$ (right).}\label{fig: comparing}
\end{figure}

Although the conclusion of the triangle recurrence may not hold when the endpoints of $e$ do not have degree $2$, there are cases when the conclusion still does hold.
For example, the graph $G_2$ in Figure~\ref{fig: comparing} cannot be constructed from $G_1$ in a way that allows us to combine the subdivision and triangle recurrences, yet we still have $|\fkD(G_2)| = 3|\fkD(G_1)|$.
Note that the conditions under discussion are local conditions; this will contrast with global conditions that we examine in Section~\ref{sec:applications}.
This leads us to ask the following.

\begin{question}
    Under what local conditions for a graph $G$ and an edge $e$ is there a ``nice'' recurrence for $\NVol(\APQ_{G \triangle e})$?
\end{question}

To close this section, we state and prove the lemmas needed to complete the proof of Theorem~\ref{thm: triangle recurrence}.

\begin{lemma}\label{lem: partial reverse 1}
	Let $G$ be any connected graph on $[N]$ for which $e = uv$ with $\deg_{G}(u) = 2$ and the neighbors of $u$ are neighbors of each other.
	If $(d_1,\dots,d_N,0) \in \fkD(G \triangle e)$, then $(d_1,\dots,d_N) - e_u \in \fkD(G)$ or $(d_1,\dots,d_N) - e_v \in \fkD(G)$.
\end{lemma}

\begin{proof}
	As usual we assume $V(G) = [N]$, $e = \{N-1,N\}$, and $\deg_{G}(N-1) = 2$.
	We will further assume that the other neighbor of $N-1$ in $G$ is $N-2$.
	
	We first show that if $d_{N-1} \geq 1$, then $(d_1,\dots,d_{N-2},d_{N-1}-1,d_N) \in \fkD(G)$.
	For notational convenience, we will write
	\[
		d' = (d_1',\dots,d_N') = (d_1,\dots,d_{N-2},d_{N-1}-1,d_N).
	\] 
	
	Consider a sum $d'_{i_1} + \cdots + d'_{i_k}$ with $1 \leq i_1 < \cdots < i_k \leq N$.
	If $i_k < N-1$, then the neighbors of each $i_j$ is the same in $D(G)$ and $D(G \triangle e)$, so
	\[
		d'_{i_1} + \cdots + d'_{i_k} = d_{i_1} + \cdots + d_{i_k} < \left| \bigcup_{j=1}^k \calN_{D(G \triangle e)}(i_j)\right| = \left| \bigcup_{j=1}^k \calN_{D(G)}(i_j)\right|.
	\]
	If $i_k = N-1$ or if both $i_k = N$ and $i_{k-1} = N-1$  then we have
	\[
		d'_{i_1} + \cdots + d'_{i_k} = d_{i_1} + \cdots + d_{i_k} - 1 < \left| \bigcup_{j=1}^k \calN_{D(G \triangle e)}(i_j)\right| - 1 = \left| \bigcup_{j=1}^k \calN_{D(G)}(i_j)\right|
	\]
	since $\overline{N+1}$ is a neighbor of $N-1$ in $D(G \triangle e)$ but not $D(G)$.
	
	Lastly, if $i_k = N$ but $i_{k-1} < N-1$, recall that we have required $N$ to be a neighbor of both $N-1$ and $N-2$.
	Thus, the neighbors of $N-1$ in $D(G \triangle e)$ are necessarily also neighbors of $N$ in $D(G \triangle e)$.
	Therefore, 
	\[\begin{aligned}
		d'_{i_1} + \cdots + d'_{i_k} &\leq d_{i_1} + \cdots + d_{i_k} + d_{N-1} - 1 \\
							&< \left|\left(\bigcup_{j=1}^k \calN_{D(G \triangle e)}(i_j)\right) \cup \calN_{D(G \triangle e)}(N-1)\right| - 1 \\
							&= \left|\bigcup_{j=1}^k \calN_{D(G \triangle e)}(i_j)\right| - 1 \\
							&= \left|\bigcup_{j=1}^k \calN_{D(G)}(i_j)\right|.
	\end{aligned}\]
	Thus, $d' \in \fkD(G)$ when $d_{N-1} \geq 1$.	
	
	If $d_{N-1} = 0$, then $d_N \geq 1$ since we may not have $d_{N-1} = d_{N} = d_{N+1} = 0$.
	We will show that, in this case, $(d_1,\dots,d_{N-1},d_N-1) \in \fkD(G)$.
	Again for notational convenience, we will write
	\[
		d'' = (d_1'',\dots,d_N'') = (d_1,\dots,d_{N-2},d_{N-1},d_N-1).
	\]
	Consider a sum $d''_{i_1} + \cdots + d''_{i_k}$ with $1 \leq i_1 < \cdots < i_k \leq N$.
	If $i_k < N-1$, then, as for $d'$,
	\[
		d''_{i_1} + \cdots + d''_{i_k} = d_{i_1} + \cdots + d_{i_k} < \left| \bigcup_{j=1}^k \calN_{D(G \triangle e)}(i_j)\right| = \left| \bigcup_{j=1}^k \calN_{D(G)}(i_j)\right|.
	\]
	If $i_k = N-1$, then since we know $d''_{N-1} = d_{N-1} = 0$, we may say
	\[
		d''_{i_1} + \cdots + d''_{i_k} = d_{i_1} + \cdots + d_{i_{k-1}} < \left| \bigcup_{j=1}^{k-1} \calN_{D(G \triangle e)}(i_j)\right| = \left| \bigcup_{j=1}^{k-1} \calN_{D(G)}(i_j)\right| \leq \left| \bigcup_{j=1}^{k} \calN_{D(G)}(i_j)\right|.
	\]
	Lastly, if $i_k = N$, then
	\[
		d''_{i_1} + \cdots + d''_{i_k} = d_{i_1} + \cdots + d_{i_k} - 1 < \left| \bigcup_{j=1}^k \calN_{D(G \triangle e)}(i_j)\right| - 1 = \left| \bigcup_{j=1}^k \calN_{D(G)}(i_j)\right|.
	\]
	Thus, $d'' \in \fkD(G)$.
	This completes the proof.
\end{proof}

\begin{lemma}\label{lem: partial reverse 2}
	Let $G$ be any connected graph on $[N]$ for which $e = uv$ with $\deg_{G}(u) = \deg_G(v) = 2$.
	If $(d_1,\dots,d_N,0) \in \fkD(G \triangle e)$, then $(d_1,\dots,d_N) - e_u \in \fkD(G)$ or $(d_1,\dots,d_N) - e_v \in \fkD(G)$.
\end{lemma}

\begin{proof}
    As usual we assume $V(G) = [N]$, $e = \{N-1,N\}$, and, this time, $\deg_{G}(N-1) = \deg_G(N) = 2$.
    If the neighbors of $N-1$ are neighbors of each other, then we are done by Lemma~\ref{lem: partial reverse 1}.
    So, we assume that the neighbors of $N-1$ are nonadjacent.
    We also may assume that the other neighbor of $N-1$ in $G$ is $N-2$ and the other neighbor of $N$ in $G$ is $N-3$.
    
    If $(d_1,\dots,d_{N-2},d_{N-1}-1,d_N) \in \fkD(G)$, then we are done.
    Otherwise, we will show $d' \in \fkD(G)$, where
    \[
    	d' = (d_1',\dots,d_N') = (d_1,\dots,d_{N-1},d_N-1)
    \]
    
    Consider a sum $d_{i_1}' + \cdots + d_{i_k}'$ with $1 \leq i_1 < \cdots < i_k \leq N$.
    First suppose $d_{N-1} = 0$.
    Since we cannot have $d_{N-1} = d_N = d_{N+1} = 0$ in a $D(G \triangle e)$-draconian sequence, it must be true that $d_N \geq 1$, so that $d'$ consists of nonnegative integers.
    
    If $i_k < N-1$, then 
    \[
    	d_{i_1}' + \cdots + d_{i_k}' = d_{i_1} + \cdots + d_{i_k} < \left|\bigcup_{j=1}^k \calN_{D(G \triangle e)}(i_j)\right| = \left|\bigcup_{j=1}^k \calN_{D(G)}(i_j)\right|
    \]
    since the neighbors of each $i_j$ are the same in $G$ and $G \triangle e$.
    If $i_k = N-1$, then $d_{i_k}' = 0$, so that 
    \[
    	d_{i_1}' + \cdots + d_{i_k}' = d_{i_1} + \cdots + d_{i_{k-1}} < \left|\bigcup_{j=1}^{k-1} \calN_{D(G \triangle e)}(i_j)\right| = \left|\bigcup_{j=1}^{k-1} \calN_{D(G)}(i_j)\right| \leq  \left|\bigcup_{j=1}^k \calN_{D(G)}(i_j)\right|.
    \]
    Lastly, if $i_k = N$, then
    \[
    	d_{i_1}' + \cdots + d_{i_k}' = d_{i_1} + \cdots + d_{i_k} - 1 < \left|\bigcup_{j=1}^k \calN_{D(G \triangle e)}(i_j)\right| -1 = \left|\bigcup_{j=1}^k \calN_{D(G)}(i_j)\right|.
    \]
    Thus, $d' \in \fkD(G)$ when $d_{N-1} = 0$.
    
    Now suppose $d_{N-1} > 0$.
    Our assumption that $(d_1,\dots,d_{N-2},d_{N-1}-1,d_N) \notin \fkD(G)$ implies $d_N > 0$ as well. 
    If $d_{N-1} = 3$, then
    \[
    	d_{N-1} + d_N < \left| \calN_{D(G \triangle e)}(N-1) \cup \calN_{D(G \triangle e)}(N)\right| = 5
    \]
    implies $d_N \leq 1$.
    We claim that this means $(d_1,\dots,d_{N-2},d_{N-1}-1,d_N) \in \fkD(G)$, which is a contradiction.
    Note that this means
    \[
    	(d_{N-1}-1,d_N) = (2,0) \text{ or } (d_{N-1}-1,d_N) = (2,1).
    \]
    
    In either situation, set $d^u =  (d^u_1,\dots,d^u_N) = (d_1,\dots,d_{N-2},d_{N-1}-1,d_N)$ and consider a sum of the form $d^u_{i_1} + \cdots + d^u_{i_k}$.
    If $i_k < N-1$, then $\overline{N+1}$ is not a neighbor of any $i_j$ in $D(G \triangle e)$, so
    \[
    	d^u_{i_1} + \cdots + d^u_{i_k} = d_{i_1} + \cdots + d_{i_k} < \left|\bigcup_{j=1}^{k} \calN_{D(G \triangle e)}(i_j)\right| = \left|\bigcup_{j=1}^{k} \calN_{D(G)}(i_j)\right|.
    \]
    If $i_k = N-1$, then $\overline{N+1}$ is a neighbor of $N-1$ in $D(G \triangle e)$, so
    \[
    	d^u_{i_1} + \cdots + d^u_{i_k} = d_{i_1} + \cdots + d_{i_k} - 1 < \left|\bigcup_{j=1}^{k} \calN_{D(G \triangle e)}(i_j)\right| - 1= \left|\bigcup_{j=1}^{k} \calN_{D(G)}(i_j)\right|.
    \]
    For the same reason, this inequality holds when $i_k = N$ and $i_{k-1} = N-1$.
    If $i_k = N$ and $i_{k-1} = N-2$, notice that the neighbors of $N-1$ are in the union of the neighbors of $N-2$ and $N$ in $D(G \triangle e)$.
    Therefore, it follows from the case in which $i_{k-1} = N-1$ that
    \[
    	\begin{aligned}
    		d^u_{i_1} + \cdots + d^u_{i_k} 
			&\leq d_{i_1} + \cdots + d_{i_k} + d_{N-1} - 1 \\
			&< \left|\bigcup_{j=1}^{k} \calN_{D(G \triangle e)}(i_j) \cup \calN_{D(G \triangle e)}(N-1)\right| - 1 \\
			&= \left|\bigcup_{j=1}^{k} \calN_{D(G)}(i_j)\right|.
	\end{aligned}
    \]
    Lastly, if $i_{k-1} < N-2$, then $\overline{N-1}$ is not a neighbor of any $i_j < N-2$ in $D(G \triangle e)$.
    Moreover, for each $i_j < N-2$, its neighbors in $D(G)$ are the same as its neighbors in $D(G \triangle e)$. 
    Putting this together with the fact that $d_N \leq 1$, we see
    \[
    	\begin{aligned}
    		d^u_{i_1} + \cdots + d^u_{i_k} &\leq d_{i_1} + \cdots + d_{i_{k-1}} + 1 \\
			&< \left|\bigcup_{j=1}^{k-1} \calN_{D(G \triangle e)}(i_j)\right| + \left|\{\overline{N-1}\}\right| \\
			&\leq \left|\bigcup_{j=1}^{k} \calN_{D(G)}(i_j)\right|.
	\end{aligned}
    \]
    Therefore, if $d_{N-1} = 3$, then $d^u \in \fkD(G)$, which is a contraction.
    
    Now suppose that $d_{N-1} = 2$.
    Analogous to before, this implies $d_N \leq 2$, leading us to the three cases
    \[
    	(d_{N-1}-1, d_N) = (1,0) \text{ or } (d_{N-1}-1, d_N) = (1,1) \text{ or } (d_{N-1}-1, d_N) = (1,2).
    \]
    If $d_N = 2$, then an argument symmetric to the one in the previous paragraph draws the same contradiction.
    If $d_N = 1$, then an argument identical to that of the previous paragraph holds.    
    Finally, if $d_N = 0$, then the desired inequalities hold since those not involving the index $N$ hold for the case of $d_N = 1$, and each inequality involving an index $N$ can be obtained from adding $d_N = 0$ to the left hand side and including $\calN_{D(G)}(N)$ in the union on the right hand side.
    Therefore, $d^u \in \fkD(G)$ whenever $d_{N-1} > 1$, which is a contradiction.

    Knowing now that $d_{N-1} = 1$, set 
    \[
    	d'' = (d_1'',\dots,d''_N) = (d_1,\dots,d_{N-2}, 1, d_N - 1)
    \]
    and consider a sum $d_{i_1}'' + \cdots + d_{i_k}''$ with $1 \leq i_1 < \cdots < i_k \leq N$.
    If $i_k = N-2$ or $i_k = N$, then the corresponding $D(G)$-inequalities hold via now-standard arguments.
    If $i_k = N-1$, then there are two subcases to consider.
    
    First suppose at least one of $\overline{N}, \overline{N-1},$ or $\overline{N-2}$ does not appear in 
    \[
    	\bigcup_{j=1}^{k-1} \calN_{D(G)}(i_j).
    \]
    Without loss of generality, assume that $\overline{N}$ does not appear.
    We can therefore say
    \[
    	d''_{i_1} + \cdots + d''_{k-1} + 1 < \left|\bigcup_{j=1}^{k-1} \calN_{D(G)}(i_j)\right| + \left| \{\overline{N}\} \right| \leq \left|\bigcup_{j=1}^k \calN_{D(G)}(i_j)\right|.
    \]
    Otherwise suppose
    \[
    	\{\overline{N-2}, \overline{N-1}, \overline{N}\} \subseteq \bigcup_{j=1}^{k-1} \calN_{D(G)}(i_j).
    \]
    Since $\deg_G(N-1) = \deg_G(N) = 2$, we know that this can only happen if $i_{k-1} = N-2$ and $i_{k-2} = N-3$.
    In particular, $i_j \neq N$ for all $j$ and
    \[
    	\calN_G(N) \subseteq \bigcup_{j=1}^{k-1} \calN_{D(G)}(i_j),
    \]
    which implies from the case $i_k = N$ that
    \[
    	\begin{aligned}
	    	d''_{i_1} + \cdots + d''_{k-1} + 1 
			&\leq d_{i_1} + \cdots + d_{k-1} + 1 + d_N - 1 \\
			&< \left|\bigcup_{j=1}^{k-1} \calN_{D(G)}(i_j) \cup \calN_{D(G)}(N)\right| \\
			&= \left| \bigcup_{j=1}^k \calN_{D(G)}(i_j)\right|.
	\end{aligned}
    \]
    This completes the proof.
\end{proof}


\subsection{Application: outerplanar graphs}\label{sec:applications}

Recall that a \emph{plane graph} is a planar graph $G$ together with a particular embedding of $G$ into the plane. 
Also recall that the \emph{weak dual} of a plane graph $G$, denoted $\wkd{G}$, is the subgraph of the dual $G^*$ induced by the vertices corresponding to bounded faces of $G$.
We denote by $E_k$ the \emph{empty graph} on $k$ vertices, that is, the disjoint union of $k$ distinct vertices.
Further, given a bounded face $F$, let $o_G(F)$ denote the number of edges of $G$ bounding both $F$ and the outer face and let $v_F$ denote the vertex of $\wkd{G}$ corresponding to $F$.
Let $\scrF(G)$ be the set of bounded faces of $G$.

\begin{definition}
	Let $G$ be a plane graph.
	The \emph{extended weak dual} of $G$, denoted $\ewd{G}$, is
	\[
		\ewd{G} = \wkd{G} \cup \left(\bigcup_{F \in \scrF(G)} v_F \vee E_{o(F)}\right)
	\] 
\end{definition}

\begin{figure}
\begin{center}
\begin{tikzpicture}
\begin{scope}[every node/.style={circle,fill=gray!50,inner sep=0pt,minimum size=2mm}]
	\node (A) at (1,1) {};
	\node (C) at (3,3) {};
	\node (D) at (3,1) {};
	\node (E) at (3,-1) {};
	\node (F) at (5,1) {} ;
\end{scope}
\begin{scope}[every node/.style={circle,fill,inner sep=0pt,minimum size=2mm}]
	\node (AA) at (2,1) {};
	\node (BB) at (3.5,1.75) {};
	\node (CC) at (3.5,0.25) {};
\end{scope}

\draw[thick, gray!50] (A) -- (C) -- (F) -- (E) -- (A);
\draw[thick, gray!50] (E) -- (C);
\draw[thick, gray!50] (D) -- (F);

\draw[thick,dashed] (AA) -- (BB) -- (CC) -- (AA);
\end{tikzpicture}
\qquad
\begin{tikzpicture}
\begin{scope}[every node/.style={circle,fill=gray!50,inner sep=0pt,minimum size=2mm}]
	\node (A) at (1,1) {};
	\node (C) at (3,3) {};
	\node (D) at (3,1) {};
	\node (E) at (3,-1) {};
	\node (F) at (5,1) {} ;
\end{scope}

\draw[thick, gray!50] (A) -- (C) -- (F) -- (E) -- (A);
\draw[thick, gray!50] (E) -- (C);
\draw[thick, gray!50] (D) -- (F);

\begin{scope}[every node/.style={circle,fill,inner sep=0pt,minimum size=2mm}]
	\node (AA) at (2,1) {};
	\node (BB) at (3.5,1.75) {};
	\node (CC) at (3.5,0.25) {};
	\node (DD) at (1.6,2.5) {};
	\node (EE) at (1.6,-0.5) {};
	\node (FF) at (4,2.5) {};
	\node (GG) at (4,-0.5) {};
\end{scope}

\draw[thick,dotted] (AA) -- (BB) -- (CC) -- (AA);
\draw[thick,dotted] (DD) -- (AA) -- (EE);
\draw[thick,dotted] (FF) -- (BB);
\draw[thick,dotted] (GG) -- (CC);
\end{tikzpicture}
\end{center}
\caption{A graph $G$ (gray) with its weak dual $\wkd{G}$ superimposed (left, dashed) and with its extended weak dual $\ewd{G}$ superimposed (right, dotted).}\label{fig: duals}
\end{figure}
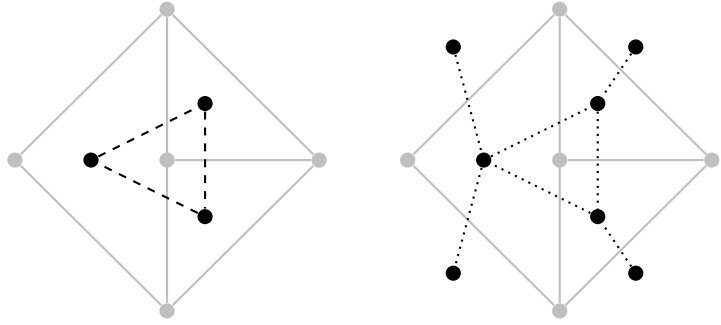

Informally, $\ewd{G}$ extends the weak dual of $G$ by including an additional edge for each edge of $G$ that bounds the outer face.
See Figure~\ref{fig: duals} for illustrations of a plane graph $G$ and its duals $\wkd{G}, \ewd{G}$.

Recall that a graph is \emph{outerplanar} if it has a planar embedding such that every vertex is incident to the outer face.
It is known \cite{outerplanarduals} that a graph is outerplanar if and only if its weak dual is a forest. 
Putting together the results of Section~\ref{sec: recurrences} we can produce a simple formula for $\NVol(\APQ_G)$ whenever $G$ can be constructed inductively by using the subdivision and triangle operations. 
The formula follows quickly from the following theorem.
	
\begin{theorem}\label{thm: ewd a tree}
	Suppose $G$ is a $2$-connected outerplane graph obtained from $C_N$ by a sequence of applications of the subdivision recurrence and the triangle recurrence.
	Then
	\[
		\NVol(\APQ_G) = 2^{1+S(G)}\prod_{F \in \scrF(G)} \deg_{\ewd{G}}(v_F),
	\]
	where
	\[
		S(G) = \sum_{F \in \scrF}(\deg_{\ewd{G}}(v_F) - 3).
	\]
\end{theorem}

\begin{proof}
	We will induct on the number of edges of $G$, which we will denote by $|E|$.
	If $G$ has $3$ or $4$ edges, then since $G$ is $2$-connected, $G = C_{|E|}$, hence $N = |E|$.
	By Corollary~\ref{cor: cycles},
	\[
		\NVol(\APQ_G) = 2^{1+(\deg_{\ewd{G}}(v_F) - 3)}\deg_{\ewd{G}}(v_F)
	\]
	as claimed.
	One may verify that this holds for $C_4$ directly as well. 
	
	Now suppose $|E| > 4$.
	If $G$ is of the form $G = G' \triangle e$ for some edge $e$ of $G'$, then let $F_0$ be the unique face of $\scrF(G) \setminus \scrF(G')$.
%
	We can then say that
	\[
		\deg_{\ewd{G'}}(v_F) = \deg_{\ewd{G}}(v_F)
	\]
	for all internal vertices $v_F$ of $\ewd{G'}$, from which it follows that $S(G) = S(G')$.
	By the triangle recurrence and the inductive assumption,
	\[\begin{aligned}
		\NVol(\APQ_G) &= 3\NVol(\APQ_{G'}) \\
			&= (\deg_{\ewd{G}}(v_{F_0}))2^{1+S(G')}\prod_{F \in \scrF(G')} \deg_{\ewd{G'}}(v_F) \\
			&= 2^{1+S(G)}\prod_{F \in \scrF(G)} \deg_{\ewd{G}}(v_F)
	\end{aligned}\]
	as desired.
	
	Suppose instead that $G$ is of the form $G = G' \col e$ for some edge $e$ of $G'$. 
	Since $G$ is outerplanar, so is $G'$, and $e$ is incident to a unique bounded face. 
	Let $B$ be the set of cut-edges of $G \setminus e$ and again let $F_0$ be the unique face of $\scrF(G) \setminus \scrF(G')$.
	The graph $H = G \setminus (B \cup \{e\})$ is the disjoint union of $k = \deg_{\ewd{G}}(v_{F_0}) - 1 - |B|$ graphs, where each component is a $2$-connected subgraph of $G$.
	Notice as well that $\scrF(H) \subseteq \scrF(G)$ with $|\scrF(H)| = |\scrF(G)| - 1$ and that
	\[
		\deg_{\ewd{H}}(v_F) = \deg_{\ewd{G}}(v_F)
	\]
	for all $F \in \scrF(H) \cap \scrF(G)$.
	By the subdivision recurrence, Corollary~\ref{cor: forest volume}, and Proposition~\ref{prop: disjoint union}, we obtain
	\begin{equation}\label{eq: outerplanar}
	\begin{aligned}
		\NVol(\APQ_G) &= \NVol(\APQ_{G' \col e}) \\
			&= 2\NVol(\APQ_{G'}) + \NVol(\APQ_{G' \setminus e}) \\
			&= 2\left(2^{1 + S(G')}\prod_{F \in \scrF(G')} \deg_{\ewd{G'}}(v_F)\right) + 2^{|B|}\NVol(\APQ_H) \\
			&= \left(2^{1+S(G)}(\deg_{\ewd{G}}(v_{F_0}) -1)\prod_{F \in \scrF(G) \setminus \{F_0\}} \deg_{\ewd{G}}(v_F)\right) \\
			& \qquad + 2^{|B|}2^{\omega}\prod_{F \in \scrF(G) \setminus \{F_0\}} \deg_{\ewd{G}}(v_F)
	\end{aligned}\end{equation}
	where 
	\[\begin{aligned}
		\omega &= k + S(G) - (\deg_{\ewd{G}}(v_{F_0})-3) \\
			&= \deg_{\ewd{G'}}(v_{F_0}) - 1 - |B| + S(G) - (\deg_{\ewd{G}}(v_{F_0})-3) \\
			&= 1 - |B| + S(G).
	\end{aligned}\]
	Simplifying the final expression in \eqref{eq: outerplanar} yields the claimed formula, completing the proof.
\end{proof}

Theorem~\ref{thm: ewd a tree} is the final piece needed to compute $\NVol(\APQ_G)$ for any outerplane graph whose $2$-connected components satisfy the conditions of Theorem~\ref{thm: ewd a tree}.

\begin{corollary}\label{cor: largest subclass}
	Let $G$ be any outerplane graph on $[N]$ such that each block with at least three vertices is obtained from $C_N$ by a sequence of applications of the subdivision recurrence and the triangle recurrence.
	Label the components of $G$ by $G_1,\dots,G_k$ and let $B_{i,1},\dots,B_{i,b_i}$ be the blocks of $G_i$.
	Then
	\begin{equation}\label{eq: outerplane}\pushQED{\qed}
		\NVol(\APQ_G) = \prod_{i=1}^k \prod_{j=1}^{b_i} 2^{1+S(B_{i,j})}\prod_{F \in \scrF(B_{i,j})} \deg_{\ewd{B_{i,j}}}(v_F). \qedhere
	\popQED\end{equation}
\end{corollary}

The graphs satisfying the conditions needed in Corollary~\ref{cor: largest subclass} form a proper, but large, class of outerplane graphs.
Experimental data suggests that the formula is, in fact, true for all outerplane graphs, but a proof eludes the authors.

\begin{conjecture}
	For \emph{any} outerplane graph $G$, Equation~\eqref{eq: outerplane} holds.
\end{conjecture}


\section{Beyond outerplanarity}

Outerplanar graphs form a large class of graphs but are far from the class of planar graphs, let alone all graphs.
For example, even though there are about $56.7 \times 10^9$ labeled outerplanar graphs on $10$ vertices, these account for only approximately $1.76\%$ of all labeled planar graphs on $10$ vertices \cite[Sequences A098000, A066537]{OEIS}. 
Because of the difficulty in computing $\NVol(\APQ_G)$ for all graphs, a natural next step would be to consider graphs that are not-quite-outerplanar. 
Toward this end, we use the following alternate characterization of outerplanar graphs.

\begin{theorem}[{\cite[Theorem~10.24]{ChartrandLesniakZhang}}]
	A graph is outerplanar if and only if contains no subdivision of $K_4$ or $K_{2,3}$ as a subgraph.
\end{theorem}

This is a direct analogue of Kuratowski's theorem, allowing one to study graphs $G$ that contain no subdivision of $K_5$ or $K_{3,3}$ but may contain a subdivision of $K_4$ or $K_{2,3}$.
In this case, a formula for $|\fkD(G)|$ remains elusive, although we do have the following partial result.
We use the notation $K^0_{M,N}$ to denote the complete bipartite graph with partite sets $[0,\dots,M-1]$ and $[M,M+N-1]$.

\begin{proposition}\label{prop: K2N}
	For all $N \geq 3$,
	\[
		\NVol(\APQ_{K_{2,N-2}}) = 2^{N-4}(N^2-N+6)-2.
	\]
\end{proposition}

\begin{proof}
	If $(c_1,\dots,c_{N}) \in \fkD(K_{2,N-2})$, then $c_1 + c_2 = k$ for some $0 \leq k \leq N-1$.
	All possible choices of $c_1,c_2$ are part of a $D(K_{2,N-2})$-draconian sequence except for $(c_1,c_2) \in \{(N-1,0),(0,N-1)\}$ since these are the only two resulting in sequences not satisfying the corresponding draconian inequalities.
	However, for the moment, we will include these in our calculations for algebraic ease.
	
	In order to satisfy the $D(K_{2,N-2})$-draconian inequalities we need the subsequence $c' = (c_3,\dots,c_{N })$ to be a weak composition of $N-1-k$ using $0$s, $1$s, and $2$s such that there is at most one $2$.
	This leads to two cases: if $c'$ contains a $2$, then there must be $N-3-k$ copies of $1$ and $k$ copies of $0$. 
	A simple counting argument gives
	\[
		(N-2)\binom{N-3}{k}
	\]
	such possibilities.
	On the other hand if $c'$ does not contain any $2$s, then there must be $N-1-k$ copies of $1$ and $k-1$ copies of $0$.
	There are $\binom{N-2}{k-1}$ such possibilities.
	Adding the values from these two cases and summing over all $k$ yields
	\[
		\sum_{k=0}^{N-1} (k+1)\left((N-2)\binom{N-3}{k}+\binom{N-2}{k-1}\right).
	\]
	The reader may verify that this simplifies to $2^{N-4}(N^2-N+6)$.
	Subtracting the two compositions where $(c_1,c_2) \in \{(N-1,0),(0,N-1)\}$ and applying Theorem~\ref{thm: translation} gives us our final formula.
\end{proof}

\begin{question}
	What is $\NVol(\APQ_{K_{M,N}})$ for arbitrary $M,N$?
\end{question}

Notice that the formula in Proposition~\ref{prop: K2N} cannot be written in the form of \eqref{eq: outerplane}.
Thus, a general formula for planar graphs will require refining the techniques of Section~\ref{sec: recurrences} or separate tools altogether.

A second important class of graphs which are planar but not outerplanar is the class of \emph{wheel graphs} $W_N = K_1 \vee C_N$.
We conjecture the following.

\begin{conjecture}\label{conj: wheels}
	For all $N\geq 3$,
	\[
		\NVol(\APQ_{W_N}) = 3^N - 2^N + 1.
	\]
\end{conjecture}

This conjecture has been verified for all $3 \leq N \leq 13$.
Wheels were examined in detail in \cite{ManyFaces} within a related, but distinct, context from $\APQ_{W_N}$.
We hope to uncover similarly rich structure in the present setting.
It may be useful to recognize that
\[
    3^N - 2^N + 1 = 2S(N+1,3) + S(N+1,2) + S(N+1,1),
\]
where $S(n,k)$ denotes the Stirling number of the second kind.

\begin{remark}
	In the time since this article was first prepared, Conjecture~\ref{conj: wheels} has been proven by Ohsugi and Tsuchiya \cite{OhsugiTsuchiya}.
\end{remark}

Finally, we give another broad class of graphs which contains all outerplanar graphs but not all planar graphs.
Strictly speaking, these graphs will allow for repeated edges, but as repeating an edge in $G$ does not affect $\APQ_G$, we need not worry about that case.

Following \cite{SeriesParallel}, first consider the directed graphs formed in the following way.
Begin with a single edge and designate one vertex the source and another vertex the sink.
This is an example of a \emph{two-terminal series-parallel graph}.
All other two-terminal series-parallel graphs are those formed by applying one of the following operations to two existing two-terminal series-parallel graphs $G$ and $H$ with sources $g$ and $h$ and sinks $g'$ and $h'$, respectively, 
\begin{enumerate}
	\item \emph{parallel composition}: produce a new graph $\mathcal{P}(G,H)$ by identifying $g$ with $h$ and $g'$ with $h'$. 
		The source of $\mathcal{P}(G,H)$ is $g \sim h$ and its sink is $g' \sim h'$.
	\item \emph{series composition}: produce a new graph $\mathcal{S}(G,H)$ by identifying $g'$ with $h$.
		The source of $\mathcal{S}(G,H)$ is $g$ and its sink is $h'$.
\end{enumerate}
A graph $G$ is a \emph{series-parallel graph} if there are two vertices $x,y$ such that, when designating $x$ as the source and $y$ as the sink, $G$ can be obtained through a sequence of applications of $\mathcal{P}(\cdot,\cdot)$ and $\mathcal{S}(\cdot,\cdot)$ when starting with a disjoint union of edges.

Series-parallel graphs are of interest in computer algorithms, as recognizing them is difficult but not intractable.
For our purposes, they are of interest because their recursive structure suggests that they may be good candidates for computing $\NVol(\APQ_G)$.
In fact, we have already seen an example of a series-parallel graph: $K_{2,N-2}$ is the parallel composition of $N-2$ copies of $P_3$, each of which is a series composition of two edges.
We ask the following question broadly, and would be interested in seeing answers to even nontrivial subclasses which are not outerplanar.

\begin{question}
	What is $\NVol(\APQ_G)$ for a series-parallel graph $G$?
\end{question}

\section{Acknowledgements}

The authors would like to thank Florian Kohl, Joakim Jakovleski, Qizhe Pan, and the anonymous referee for their detailed feedback.
Their comments greatly improved the quality of this article.


\end{document}